


\input amstex
\expandafter\ifx\csname mathdefs.tex\endcsname\relax
  \expandafter\gdef\csname mathdefs.tex\endcsname{}
\else \message{Hey!  Apparently you were trying to
  \string twice.   This does not make sense.} 
\errmessage{Please edit your file (probably \jobname.tex) and remove
any duplicate ``\string\input'' lines} \fi




\catcode`\X=12\catcode`\@=11

\def\n@wcount{\alloc@0\count\countdef\insc@unt}
\def\n@wwrite{\alloc@7\write\chardef\sixt@@n}
\def\n@wread{\alloc@6\read\chardef\sixt@@n}
\def\r@s@t{\relax}\def\v@idline{\par}\def\@mputate#1/{#1}
\def\l@c@l#1X{\firstpart.#1}\def\gl@b@l#1X{#1}\def\t@d@l#1X{{}}

\def\crossrefs#1{\ifx\all#1\let\tr@ce=\all\else\def\tr@ce{#1,}\fi
   \n@wwrite\cit@tionsout\openout\cit@tionsout=\jobname.cit 
   \write\cit@tionsout{\tr@ce}\expandafter\setfl@gs\tr@ce,}
\def\setfl@gs#1,{\def\@{#1}\ifx\@\empty\let\next=\relax
   \else\let\next=\setfl@gs\expandafter\xdef
   \csname#1tr@cetrue\endcsname{}\fi\next}
\def\m@ketag#1#2{\expandafter\n@wcount\csname#2tagno\endcsname
     \csname#2tagno\endcsname=0\let\tail=\all\xdef\all{\tail#2,}
   \ifx#1\l@c@l\let\tail=\r@s@t\xdef\r@s@t{\csname#2tagno\endcsname=0\tail}\fi
   \expandafter\gdef\csname#2cite\endcsname##1{\expandafter
     \ifx\csname#2tag##1\endcsname\relax?\else\csname#2tag##1\endcsname\fi
     \expandafter\ifx\csname#2tr@cetrue\endcsname\relax\else
     \write\cit@tionsout{#2tag ##1 cited on page \folio.}\fi}
   \expandafter\gdef\csname#2page\endcsname##1{\expandafter
     \ifx\csname#2page##1\endcsname\relax?\else\csname#2page##1\endcsname\fi
     \expandafter\ifx\csname#2tr@cetrue\endcsname\relax\else
     \write\cit@tionsout{#2tag ##1 cited on page \folio.}\fi}
   \expandafter\gdef\csname#2tag\endcsname##1{\expandafter
      \ifx\csname#2check##1\endcsname\relax
      \expandafter\xdef\csname#2check##1\endcsname{}%
      \else\immediate\write16{Warning: #2tag ##1 used more than once.}\fi
      \multit@g{#1}{#2}##1/X%
      \write\t@gsout{#2tag ##1 assigned number \csname#2tag##1\endcsname\space
      on page \number\count0.}%
   \csname#2tag##1\endcsname}}
\def\multit@g#1#2#3/#4X{\def\t@mp{#4}\ifx\t@mp\empty%
      \global\advance\csname#2tagno\endcsname by 1 
      \expandafter\xdef\csname#2tag#3\endcsname
      {#1\number\csname#2tagno\endcsnameX}%
   \else\expandafter\ifx\csname#2last#3\endcsname\relax
      \expandafter\n@wcount\csname#2last#3\endcsname
      \global\advance\csname#2tagno\endcsname by 1 
      \expandafter\xdef\csname#2tag#3\endcsname
      {#1\number\csname#2tagno\endcsnameX}
      \write\t@gsout{#2tag #3 assigned number \csname#2tag#3\endcsname\space
      on page \number\count0.}\fi
   \global\advance\csname#2last#3\endcsname by 1
   \def\t@mp{\expandafter\xdef\csname#2tag#3/}%
   \expandafter\t@mp\@mputate#4\endcsname
   {\csname#2tag#3\endcsname\lastpart{\csname#2last#3\endcsname}}\fi}
\def\t@gs#1{\def\all{}\m@ketag#1e\m@ketag#1s\m@ketag\t@d@l p
   \m@ketag\gl@b@l r \n@wread\t@gsin
   \openin\t@gsin=\jobname.tgs \re@der \closein\t@gsin
   \n@wwrite\t@gsout\openout\t@gsout=\jobname.tgs }
\outer\def\localtags{\t@gs\l@c@l}
\outer\def\globaltags{\t@gs\gl@b@l}
\outer\def\newlocaltag#1{\m@ketag\l@c@l{#1}}
\outer\def\newglobaltag#1{\m@ketag\gl@b@l{#1}}

\newif\ifpr@ 
\def\m@kecs #1tag #2 assigned number #3 on page #4.%
   {\expandafter\gdef\csname#1tag#2\endcsname{#3}
   \expandafter\gdef\csname#1page#2\endcsname{#4}
   \ifpr@\expandafter\xdef\csname#1check#2\endcsname{}\fi}
\def\re@der{\ifeof\t@gsin\let\next=\relax\else
   \read\t@gsin to\t@gline\ifx\t@gline\v@idline\else
   \expandafter\m@kecs \t@gline\fi\let \next=\re@der\fi\next}
\def\pretags#1{\pr@true\pret@gs#1,,}
\def\pret@gs#1,{\def\@{#1}\ifx\@\empty\let\n@xtfile=\relax
   \else\let\n@xtfile=\pret@gs \openin\t@gsin=#1.tgs \message{#1} \re@der 
   \closein\t@gsin\fi \n@xtfile}

\newcount\sectno\sectno=0\newcount\subsectno\subsectno=0
\newif\ifultr@local \def\ultralocal{\ultr@localtrue}
\def\firstpart{\number\sectno}
\def\lastpart#1{\ifcase#1 \or a\or b\or c\or d\or e\or f\or g\or h\or 
   i\or k\or l\or m\or n\or o\or p\or q\or r\or s\or t\or u\or v\or w\or 
   x\or y\or z \fi}

\def\resetall{\global\advance\sectno by 1\subsectno=0
   \gdef\firstpart{\number\sectno}\r@s@t}
\def\resetsub{\global\advance\subsectno by 1
   \gdef\firstpart{\number\sectno.\number\subsectno}\r@s@t}
\def\newsection#1\par{\resetall\vskip0pt plus.3\vsize\penalty-250
   \vskip0pt plus-.3\vsize\bigskip\bigskip
   \message{#1}\leftline{\bf#1}\nobreak\bigskip}
\def\subsection#1\par{\ifultr@local\resetsub\fi
   \vskip0pt plus.2\vsize\penalty-250\vskip0pt plus-.2\vsize
   \bigskip\smallskip\message{#1}\leftline{\bf#1}\nobreak\medskip}

\def\t@gsoff#1,{\def\@{#1}\ifx\@\empty\let\next=\relax\else\let\next=\t@gsoff
   \def\@@{p}\ifx\@\@@\else
   \expandafter\gdef\csname#1cite\endcsname##1{\zeigen{##1}}
   \expandafter\gdef\csname#1page\endcsname##1{?}
   \expandafter\gdef\csname#1tag\endcsname##1{\zeigen{##1}}\fi\fi\next}
\def\verbatimtags{\ifx\all\relax\else\expandafter\t@gsoff\all,\fi}
\def\zeigen#1{\hbox{$\langle$}#1\hbox{$\rangle$}}

\def\(#1){\edef\dot@g{\ifmmode\ifinner(\hbox{\noexpand\etag{#1}})
   \else\noexpand\eqno(\hbox{\noexpand\etag{#1}})\fi
   \else(\noexpand\ecite{#1})\fi}\dot@g}

\newif\ifbr@ck
\def\eat#1{}
\def\[#1]{\br@cktrue[\br@cket#1'X]}
\def\br@cket#1'#2X{\def\temp{#2}\ifx\temp\empty\let\next\eat
   \else\let\next\br@cket\fi
   \ifbr@ck\br@ckfalse\br@ck@t#1,X\else\br@cktrue#1\fi\next#2X}
\def\br@ck@t#1,#2X{\def\temp{#2}\ifx\temp\empty\let\neext\eat
   \else\let\neext\br@ck@t\def\temp{,}\fi
   \def\teemp{#1}\ifx\teemp\empty\else\rcite{#1}\fi\temp\neext#2X}
\def\resetbr@cket{\gdef\[##1]{[\rtag{##1}]}}
\def\references{\resetbr@cket\newsection References\par}

\newtoks\symb@ls\newtoks\s@mb@ls\newtoks\p@gelist\n@wcount\ftn@mber
    \ftn@mber=1\newif\ifftn@mbers\ftn@mbersfalse\newif\ifbyp@ge\byp@gefalse
\def\defm@rk{\ifftn@mbers\n@mberm@rk\else\symb@lm@rk\fi}
\def\n@mberm@rk{\xdef\m@rk{{\the\ftn@mber}}%
    \global\advance\ftn@mber by 1 }
\def\rot@te#1{\let\temp=#1\global#1=\expandafter\r@t@te\the\temp,X}
\def\r@t@te#1,#2X{{#2#1}\xdef\m@rk{{#1}}}
\def\b@@st#1{{$^{#1}$}}\def\str@p#1{#1}
\def\symb@lm@rk{\ifbyp@ge\rot@te\p@gelist\ifnum\expandafter\str@p\m@rk=1 
    \s@mb@ls=\symb@ls\fi\write\f@nsout{\number\count0}\fi \rot@te\s@mb@ls}
\def\byp@ge{\byp@getrue\n@wwrite\f@nsin\openin\f@nsin=\jobname.fns 
    \n@wcount\currentp@ge\currentp@ge=0\p@gelist={0}
    \re@dfns\closein\f@nsin\rot@te\p@gelist
    \n@wread\f@nsout\openout\f@nsout=\jobname.fns }
\def\m@kelist#1X#2{{#1,#2}}
\def\re@dfns{\ifeof\f@nsin\let\next=\relax\else\read\f@nsin to \f@nline
    \ifx\f@nline\v@idline\else\let\t@mplist=\p@gelist
    \ifnum\currentp@ge=\f@nline
    \global\p@gelist=\expandafter\m@kelist\the\t@mplistX0
    \else\currentp@ge=\f@nline
    \global\p@gelist=\expandafter\m@kelist\the\t@mplistX1\fi\fi
    \let\next=\re@dfns\fi\next}
\def\symbols#1{\symb@ls={#1}\s@mb@ls=\symb@ls} 
\def\bigsymbol{\textstyle}
\symbols{\bigsymbol\ast,\dagger,\ddagger,\sharp,\flat,\natural,\star}
\def\ftnumbers{\ftn@mberstrue} \def\ftsymbols{\ftn@mbersfalse}
\def\paginal{\byp@ge} \def\resetftnumbers{\ftn@mber=1}
\def\ftnote#1{\defm@rk\expandafter\expandafter\expandafter\footnote
    \expandafter\b@@st\m@rk{#1}}

\long\def\jump#1\endjump{}
\def\ssum{\mathop{\lower .1em\hbox{$\textstyle\Sigma$}}\nolimits}

\def\qed{\nobreak\kern 1em \vrule height .5em width .5em depth 0em}
\def\newneq{\hbox{\rlap{\hbox to 1\wd9{\hss$=$\hss}}\raise .1em 
   \hbox to 1\wd9{\hss$\scriptscriptstyle/$\hss}}}
\def\subsetne{\setbox9 = \hbox{$\subset$}\mathrel{\hbox{\rlap
   {\lower .4em \newneq}\raise .13em \hbox{$\subset$}}}}
\def\supsetne{\setbox9 = \hbox{$\subset$}\mathrel{\hbox{\rlap
   {\lower .4em \newneq}\raise .13em \hbox{$\supset$}}}}

\def\vbar{\mathchoice{\vrule height6.3ptdepth-.5ptwidth.8pt\kern-.8pt}
   {\vrule height6.3ptdepth-.5ptwidth.8pt\kern-.8pt}
   {\vrule height4.1ptdepth-.35ptwidth.6pt\kern-.6pt}
   {\vrule height3.1ptdepth-.25ptwidth.5pt\kern-.5pt}}
\def\f@dge{\mathchoice{}{}{\mkern.5mu}{\mkern.8mu}}
\def\b@c#1#2{{\rm \mkern#2mu\vbar\mkern-#2mu#1}}
\def\b@b#1{{\rm I\mkern-3.5mu #1}}
\def\b@a#1#2{{\rm #1\mkern-#2mu\f@dge #1}}
\def\bb#1{{\count4=`#1 \advance\count4by-64 \ifcase\count4\or\b@a A{11.5}\or
   \b@b B\or\b@c C{5}\or\b@b D\or\b@b E\or\b@b F \or\b@c G{5}\or\b@b H\or
   \b@b I\or\b@c J{3}\or\b@b K\or\b@b L \or\b@b M\or\b@b N\or\b@c O{5} \or
   \b@b P\or\b@c Q{5}\or\b@b R\or\b@a S{8}\or\b@a T{10.5}\or\b@c U{5}\or
   \b@a V{12}\or\b@a W{16.5}\or\b@a X{11}\or\b@a Y{11.7}\or\b@a Z{7.5}\fi}}

\catcode`\X=11 \catcode`\@=12

\expandafter\ifx\csname citeadd.tex\endcsname\relax
\expandafter\gdef\csname citeadd.tex\endcsname{}
\else \message{Hey!  Apparently you were trying to
\string twice.   This does not make sense.} 
\errmessage{Please edit your file (probably \jobname.tex) and remove
any duplicate ``\string\input'' lines} \fi

\def\sciteu{\sciteerror{undefined}}

\def\sciteerror#1#2{{\mathortextbf{\scite{#2}}}\complainaboutcitation{#1}{#2}}
\def\mathortextbf#1{\hbox{\bf #1}}
\def\complainaboutcitation#1#2{%
\vadjust{\line{\llap{---$\!\!>$ }\qquad scite$\{$#2$\}$ #1\hfil}}}

\sectno=-1   
\localtags
\NoBlackBoxes
\define\mr{\medskip\roster}
\define\sn{\smallskip\noindent}
\define\mn{\medskip\noindent}
\define\bn{\bigskip\noindent}
\define\ub{\underbar}
\define\wilog{\text{without loss of generality}}
\define\ermn{\endroster\medskip\noindent}
\define\dbca{\dsize\bigcap}
\define\dbcu{\dsize\bigcup}
\define \nl{\newline}
\ifx\shlhetal\undefinedcontrolsequence\let\shlhetal\relax\fi
\documentstyle {amsppt}
\topmatter
\title {WAS SIERPI\'NSKI RIGHT IV? \\
Sh546} \endtitle
\author {Saharon Shelah \thanks {\null\newline
Typed 5/92 - ($2^{\aleph_0},k^2_2$-Mahlo,$\lambda \rightarrow [\aleph_2]
^2_3$; some on models) \null\newline
I thank Alice Leonhardt for the beautiful typing. \null\newline
Latest Revision 97/Dec/4 - Revised and Expanded 5/94 based on lectures,
\null\newline
Summer '94, Jerusalem. \null\newline
Partially supported by the basic research fund, Israeli
Academy of Sciences.} \endthanks} \endauthor
\affil {Institute of Mathematics \\
The Hebrew University \\
Jerusalem, Israel
\medskip
Rutgers University \\
Department of Mathematics \\
New Brunswick, NJ  USA} \endaffil
\bn
\abstract {We prove for any $\mu  = \mu^{<\mu} < \theta < \lambda,\lambda$
large enough (just strongly inaccessible Mahlo) the consistency of 
$2^\mu = \lambda \rightarrow [\theta]^2_3$ and even $2^\mu = \lambda
\rightarrow [\theta]^2_{\sigma,2}$ for $\sigma < \mu$.  The new point is
that possibly $\theta > \mu^+$.} \endabstract
\endtopmatter
\document  

\expandafter\ifx\csname alice2jlem.tex\endcsname\relax
  \expandafter\gdef\csname alice2jlem.tex\endcsname{}
\else \message{Hey!  Apparently you were trying to
\string  twice.   This does not make sense.}
\errmessage{Please edit your file (probably \jobname.tex) and remove
any duplicate ``\string\input'' lines} \fi

\expandafter\ifx\csname bib4plain.tex\endcsname\relax
  \expandafter\gdef\csname bib4plain.tex\endcsname{}
\else \message{Hey!  Apparently you were trying to \string twice.   This does not make sense.}
\errmessage{Please edit your file (probably \jobname.tex) and remove
any duplicate ``\string\input'' lines} \fi

\def\renewcommand{\newcommand}	       
\edef\cite{\the\catcode`@}%
\catcode`@ = 11
\let\@oldatcatcode = \cite
\chardef\@letter = 11
\chardef\@other = 12
%
%
%
%
\def\@innerdef#1#2{\edef#1{\expandafter\noexpand\csname #2\endcsname}}%
%
%
\@innerdef\@innernewcount{newcount}%
\@innerdef\@innernewdimen{newdimen}%
\@innerdef\@innernewif{newif}%
\@innerdef\@innernewwrite{newwrite}%
%
%
%
\def\@gobble#1{}%
%
%
%
\ifx\inputlineno\@undefined
   \let\@linenumber = \empty 
\else
   \def\@linenumber{\the\inputlineno:\space}%
\fi
%
%
%
\def\@futurenonspacelet#1{\def\cs{#1}%
   \afterassignment\@stepone\let\@nexttoken=
}%
\begingroup 
\def\\{\global\let\@stoken= }%
\\ 
\endgroup
\def\@stepone{\expandafter\futurelet\cs\@steptwo}%
\def\@steptwo{\expandafter\ifx\cs\@stoken\let\@@next=\@stepthree
   \else\let\@@next=\@nexttoken\fi \@@next}%
\def\@stepthree{\afterassignment\@stepone\let\@@next= }%
%
%
%
\def\@getoptionalarg#1{%
   \let\@optionaltemp = #1%
   \let\@optionalnext = \relax
   \@futurenonspacelet\@optionalnext\@bracketcheck
}%
%
%
\def\@bracketcheck{%
   \ifx [\@optionalnext
      \expandafter\@@getoptionalarg
   \else
      \let\@optionalarg = \empty
      \expandafter\@optionaltemp
   \fi
}%
\def\@@getoptionalarg[#1]{%
   \def\@optionalarg{#1}%
   \@optionaltemp
}%
%
%
%
\def\@nnil{\@nil}%
\def\@fornoop#1\@@#2#3{}%
\def\@for#1:=#2\do#3{%
   \edef\@fortmp{#2}%
   \ifx\@fortmp\empty \else
      \expandafter\@forloop#2,\@nil,\@nil\@@#1{#3}%
   \fi
}%
\def\@forloop#1,#2,#3\@@#4#5{\def#4{#1}\ifx #4\@nnil \else
       #5\def#4{#2}\ifx #4\@nnil \else#5\@iforloop #3\@@#4{#5}\fi\fi
}%
\def\@iforloop#1,#2\@@#3#4{\def#3{#1}\ifx #3\@nnil
       \let\@nextwhile=\@fornoop \else
      #4\relax\let\@nextwhile=\@iforloop\fi\@nextwhile#2\@@#3{#4}%
}%
%
%
%
\@innernewif\if@fileexists
\def\@testfileexistence{\@getoptionalarg\@finishtestfileexistence}%
\def\@finishtestfileexistence#1{%
   \begingroup
      \def\extension{#1}%
      \immediate\openin0 =
         \ifx\@optionalarg\empty\jobname\else\@optionalarg\fi
         \ifx\extension\empty \else .#1\fi
         \space
      \ifeof 0
         \global\@fileexistsfalse
      \else
         \global\@fileexiststrue
      \fi
      \immediate\closein0
   \endgroup
}%
%
%
%
%
\def\bibliographystyle#1{%
   \@readauxfile
   \@writeaux{\string\bibstyle{#1}}%
}%
\let\bibstyle = \@gobble
%
%
\let\bblfilebasename = \jobname
\def\bibliography#1{%
   \@readauxfile
   \@writeaux{\string\bibdata{#1}}%
   \@testfileexistence[\bblfilebasename]{bbl}%
   \if@fileexists
      \nobreak
      \@readbblfile
   \fi
}%
\let\bibdata = \@gobble
%
%
\def\nocite#1{%
   \@readauxfile
   \@writeaux{\string\citation{#1}}%
}%
\@innernewif\if@notfirstcitation
%
%
\def\cite{\@getoptionalarg\@cite}%
%
%
\def\@cite#1{%
   \let\@citenotetext = \@optionalarg
   \printcitestart
   \nocite{#1}%
   \@notfirstcitationfalse
   \@for \@citation :=#1\do
   {%
      \expandafter\@onecitation\@citation\@@
   }%
   \ifx\empty\@citenotetext\else
      \printcitenote{\@citenotetext}%
   \fi
   \printcitefinish
}%
\def\@onecitation#1\@@{%
   \if@notfirstcitation
      \printbetweencitations
   \fi
   \expandafter \ifx \csname\@citelabel{#1}\endcsname \relax
      \if@citewarning
         \message{\@linenumber Undefined citation `#1'.}%
      \fi
      \expandafter\gdef\csname\@citelabel{#1}\endcsname{%
\strut
\vadjust{\vskip-\dp\strutbox
\vbox to 0pt{\vss\parindent0cm \leftskip=\hsize 
\advance\leftskip3mm
\advance\hsize 4cm\strut\openup-4pt 
\rightskip 0cm plus 1cm minus 0.5cm ?  #1 ?\strut}}
         {\tt
            \escapechar = -1
            \nobreak\hskip0pt
            \expandafter\string\csname#1\endcsname
            \nobreak\hskip0pt
         }%
      }%
   \fi
   \csname\@citelabel{#1}\endcsname
   \@notfirstcitationtrue
}%
%
%
\def\@citelabel#1{b@#1}%
%
%
\def\@citedef#1#2{\expandafter\gdef\csname\@citelabel{#1}\endcsname{#2}}%
%
%
%
\def\@readbblfile{%
   \ifx\@itemnum\@undefined
      \@innernewcount\@itemnum
   \fi
   \begingroup
      \def\begin##1##2{%
         \setbox0 = \hbox{\biblabelcontents{##2}}%
         \biblabelwidth = \wd0
      }%
      \def\end##1{}
      %
      %
      \@itemnum = 0
      \def\bibitem{\@getoptionalarg\@bibitem}%
      \def\@bibitem{%
         \ifx\@optionalarg\empty
            \expandafter\@numberedbibitem
         \else
            \expandafter\@alphabibitem
         \fi
      }%
      \def\@alphabibitem##1{%
         \expandafter \xdef\csname\@citelabel{##1}\endcsname {\@optionalarg}%
         \ifx\biblabelprecontents\@undefined
            \let\biblabelprecontents = \relax
         \fi
         \ifx\biblabelpostcontents\@undefined
            \let\biblabelpostcontents = \hss
         \fi
         \@finishbibitem{##1}%
      }%
      \def\@numberedbibitem##1{%
         \advance\@itemnum by 1
         \expandafter \xdef\csname\@citelabel{##1}\endcsname{\number\@itemnum}%
         \ifx\biblabelprecontents\@undefined
            \let\biblabelprecontents = \hss
         \fi
         \ifx\biblabelpostcontents\@undefined
            \let\biblabelpostcontents = \relax
         \fi
         \@finishbibitem{##1}%
      }%
      \def\@finishbibitem##1{%
         \biblabelprint{\csname\@citelabel{##1}\endcsname}%
         \@writeaux{\string\@citedef{##1}{\csname\@citelabel{##1}\endcsname}}%
         \ignorespaces
      }%
      %
      %
      \let\em = \bblem
      \let\newblock = \bblnewblock
      \let\sc = \bblsc
      \frenchspacing
      \clubpenalty = 4000 \widowpenalty = 4000
      \tolerance = 10000 \hfuzz = .5pt
      \everypar = {\hangindent = \biblabelwidth
                      \advance\hangindent by \biblabelextraspace}%
      \bblrm
      \parskip = 1.5ex plus .5ex minus .5ex
      \biblabelextraspace = .5em
      \bblhook
      \input \bblfilebasename.bbl
   \endgroup
}%
%
%
\@innernewdimen\biblabelwidth
\@innernewdimen\biblabelextraspace
%
%
%
\def\biblabelprint#1{%
   \noindent
   \hbox to \biblabelwidth{%
      \biblabelprecontents
      \biblabelcontents{#1}%
      \biblabelpostcontents
   }%
   \kern\biblabelextraspace
}%
%
%
%
\def\biblabelcontents#1{{\bblrm [#1]}}%
%
%
\def\bblrm{\rm}%
%
%
\def\bblem{\it}%
%
%
\def\bblsc{\ifx\@scfont\@undefined
              \font\@scfont = cmcsc10
           \fi
           \@scfont
}%
%
%
\def\bblnewblock{\hskip .11em plus .33em minus .07em }%
%
%
\let\bblhook = \empty
%
%
%
\def\printcitestart{[}
\def\printcitefinish{]}
\def\printbetweencitations{, }
\def\printcitenote#1{, #1}
%
%
%
\let\citation = \@gobble
%
%
%
\@innernewcount\@numparams
%
%
\def\newcommand#1{%
   \def\@commandname{#1}%
   \@getoptionalarg\@continuenewcommand
}%
%
%
\def\@continuenewcommand{%
   \@numparams = \ifx\@optionalarg\empty 0\else\@optionalarg \fi \relax
   \@newcommand
}%
%
%
\def\@newcommand#1{%
   \def\@startdef{\expandafter\edef\@commandname}%
   \ifnum\@numparams=0
      \let\@paramdef = \empty
   \else
      \ifnum\@numparams>9
         \errmessage{\the\@numparams\space is too many parameters}%
      \else
         \ifnum\@numparams<0
            \errmessage{\the\@numparams\space is too few parameters}%
         \else
            \edef\@paramdef{%
               \ifcase\@numparams
                  \empty  No arguments.
               \or ####1%
               \or ####1####2%
               \or ####1####2####3%
               \or ####1####2####3####4%
               \or ####1####2####3####4####5%
               \or ####1####2####3####4####5####6%
               \or ####1####2####3####4####5####6####7%
               \or ####1####2####3####4####5####6####7####8%
               \or ####1####2####3####4####5####6####7####8####9%
               \fi
            }%
         \fi
      \fi
   \fi
   \expandafter\@startdef\@paramdef{#1}%
}%
%
%
%
%
\def\@readauxfile{%
   \if@auxfiledone \else 
      \global\@auxfiledonetrue
      \@testfileexistence{aux}%
      \if@fileexists
         \begingroup
            \endlinechar = -1
            \catcode`@ = 11
            \input \jobname.aux
         \endgroup
      \else
         \message{\@undefinedmessage}%
         \global\@citewarningfalse
      \fi
      \immediate\openout\@auxfile = \jobname.aux
   \fi
}%
%
%
\newif\if@auxfiledone
\ifx\noauxfile\@undefined \else \@auxfiledonetrue\fi
%
%
%
%
\@innernewwrite\@auxfile
\def\@writeaux#1{\ifx\noauxfile\@undefined \write\@auxfile{#1}\fi}%
%
%
%
\ifx\@undefinedmessage\@undefined
   \def\@undefinedmessage{No .aux file; I won't give you warnings about
                          undefined citations.}%
\fi
%
%
\@innernewif\if@citewarning
\ifx\noauxfile\@undefined \@citewarningtrue\fi
%
%
%
\catcode`@ = \@oldatcatcode


\def\widestnumber#1#2{}

\def\rm{\fam0 \tenrm}

\def\fakesubhead#1\endsubhead{\bigskip\noindent{\bf#1}\par}


%
%
%

%

\font\textrsfs=rsfs10
\font\scriptrsfs=rsfs7
\font\scriptscriptrsfs=rsfs5

\newfam\rsfsfam
\textfont\rsfsfam=\textrsfs
\scriptfont\rsfsfam=\scriptrsfs
\scriptscriptfont\rsfsfam=\scriptscriptrsfs

\edef\oldcatcodeofat{\the\catcode`\@}
\catcode`\@11

\def\Cal@@#1{\noaccents@ \fam \rsfsfam #1}

\catcode`\@\oldcatcodeofat

\newpage

\head {\S0 Introduction} \endhead  \resetall 
\bigskip

An important theme is modern set theory is to prove the consistency
of ``small cardinals" having ``a large cardinal property".  Probably
the dominant interpretation concerns large ideals (with reflection
properties or connected to generic embedding).  But here we deal with
another important interpretation:
partition properties.  We continue here \cite[\S2]{Sh:276}, \cite{Sh:288}, 
\cite{Sh:289}, \cite{Sh:473}, \cite{Sh:481} but generally do not rely on 
them except in the end (of the proof of \scite{1.12}) when it becomes like the
proof of \cite[\S2]{Sh:276}.  
This work is continued in Rabus and Shelah \newline
\cite{RbSh:585}.
\bigskip

We thank the participants in a logic seminar in The Hebrew University,
Spring '94, and Mariusz Rabus for their comments.
\bigskip

\centerline {Preliminaries} 
\bigskip
\noindent
\subhead {\stag{0.A}} \endsubhead  $Let  <^\ast_\chi$ be a well ordering of
\newline
${\Cal H}(\chi ) = \{ x$ : the 
transitive closure of $x$ has cardinality  $<\chi \}$  
agreeing with the usual well ordering of the ordinals, \newline
$P$  (and  $Q,R$)  
will denote forcing notions, i.e. quasi orders with a minimal element  
$\emptyset  = \emptyset_P$. 

A forcing notion $P$ is $\lambda$-closed or $\lambda$-complete if every 
increasing sequence of members of $P$, of length less than $\lambda$, has an 
upper bound.
\bigskip

\subhead {\stag{0.B}} \endsubhead  If  $P \in  {\Cal H}(\chi)$, 
then for a sequence $\bar p = \langle p_i:i < \gamma \rangle$ of members of  
$P$ (not necessarily increasing) let $\alpha = \alpha_{\bar p} =: 
\text{ sup} \bigl\{ j:\{p_i:i < j\}$ has an upper bound in $P \bigr\}$ 
and define the canonical upper bound of $\bar p$, denoted by $\& \bar p$  
as follows: 
\medskip
\roster
\item "{(a)}"  the least upper bound of $\{p_i:i < \alpha_{\bar p}\}$ 
in $P$ if there exists such an element
\sn
\item "{(b)}"  the $<^\ast_\chi$-first upper bound of  $\bar p$  if
(a) can't be applied but there is an upper bound of $\{p_i:i < 
\alpha_{\bar p}\}$,
\sn
\item "{(c)}"  $p_0$ if (a), (b) fail,  $\gamma > 0$,
\sn
\item "{(d)}"  $\emptyset_P$ if  $\gamma = 0$.
\endroster
\medskip 

\noindent
Let $p_0 \and p_1$ be the canonical upper bound of $\langle p_\ell:\ell < 2 
\rangle$. \newline
Take $[a]^\kappa  = \{ b \subseteq a:|b| = \kappa \}$ and $[a]^{<\kappa} 
= \dsize \bigcup_{\theta < \kappa} [a]^\theta$.
\bigskip

\subhead {\stag{0.C}} \endsubhead  For sets of ordinals, $A$ and $B$,
define OP$_{B,A}$ as the maximal order preserving 1-to-1 function between 
initial segments of $A$ and $B$,  
i.e., it is the function with domain  $\{\alpha \in A:\text{otp}
(\alpha \cap A) < \text{otp}(B)\}$ and $OP_{B,A}(\alpha) = \beta$ if and 
only if $\alpha \in A,\beta \in B$ and otp$(\alpha \cap A) = 
\text{ otp}(\beta \cap B)$.
\mn
If $A,B$ are sets of ordinals, let $A \triangleleft B$ mean $A$ is a proper
initial segment of $B$.  If $\eta,\nu$ are sequences let $\eta \triangleleft 
\nu$ mean $\nu$ is an initial segment of $\nu$.  If we write 
$\trianglelefteq$ (rather than $\triangleleft$) we allow equality.
\mn
Let $S^\lambda_\kappa = \{\delta < \lambda:\text{cf}(\delta) = \kappa\}$.
\newpage

\definition{Definition \stag{0.1}}  $\lambda \rightarrow [\alpha]^n_\theta$
holds provided that whenever $F$ is a function from $[\lambda]^n$ to $\theta$,
\ub{then} there is $A \subseteq \lambda$ of order type $\alpha$ and
$t < \theta$ such that $[w \in [A]^n \Rightarrow F(w) \ne t]$.
\enddefinition
\bigskip

\definition{Definition \stag{0.2}}  
$\lambda \rightarrow [\alpha]^n_{\kappa,\theta}$ if for every function $F$ 
from $[\lambda]^n$ to $\kappa$ there is $A \subseteq \lambda$ of order type  
$\alpha$ such that $\{F(w):w \in [A]^n\}$ has power $\le \theta$.  If we write
``$< \theta$" instead of $\theta$ we mean that the set above has cardinality
$< \theta$.
\enddefinition
\bigskip

\definition{Definition \stag{0.3}}  A forcing notion  $P$  satisfies the
Knaster condition (has property $K$) if for any 
$\{p_i:i < \omega_1\} \subseteq P$  there is an uncountable 
$A \subseteq \omega_1$ such that the conditions $p_i$ and $p_j$ are 
compatible whenever $i$, $j \in A$.
\enddefinition
\bigskip

\noindent
What problems do \cite{Sh:276}, \cite{Sh:288}, \cite{Sh:289}, \cite{Sh:473} 
and \cite{Sh:481} raise?  The most important ``minimal open", as suggested 
in \cite{Sh:481} were:
\subhead {A Question} \endsubhead  (1)  Can we get e.g. CON$(2^{\aleph_0}
\rightarrow [\aleph_2]^2_3)$ (generally raise $\mu^+$ in part (3) below
to higher cardinals).  We solve it here. \newline
(2)  Can we get CON$(\aleph_\omega > 2^{\aleph_0} \rightarrow 
[\aleph_1]^2_3)$ (the exact $\aleph_n$ seems to me less exciting). \newline
(3)  Can we get e.g. CON$(2^\mu > \lambda \rightarrow [\mu^+]^2_3$)?
\medskip

\noindent
Also
\subhead {B Question} \endsubhead  (1) Can we get the continuity on a
non-meagre set for functions $f:{}^\kappa 2 \rightarrow {}^\kappa 2$? 
(Solved in \cite{Sh:473}.)  \newline
(2) What can we say on continuity of $2$-place functions (dealt with in
Rabus Shelah \cite{RbSh:585})? \newline
(3) What about $n$-place functions?
(continuing in this respect \cite{Sh:288} probably just combine
\cite{RbSh:585} with)
\medskip

\subhead {C Question} \endsubhead  (1) \cite{Sh:481} for 
$\mu > \aleph_0$. \newline
(2) Can we get e.g. CON($2^{\aleph_0} \ge \aleph_2$, and if $P$ is 
$2^{\aleph_0}$-c.c., $Q$ is $\aleph_2$-c.c., then $P \times Q$ is 
$2^{\aleph_0}$-c.c.). \nl
(3) Can we get e.g. CON$(2^{\aleph_0} > \lambda > \aleph_0$, and if $P$ is
$\lambda$-c.c., $Q$ is $\aleph_1$-c.c. then $P \times Q$ is $\lambda$-c.c.);
more general is CON($\mu = \mu^{< \mu} > \aleph_0 +$ if $P$ is $2^\mu$-c.c. 
$Q$ is $\mu^+$-c.c. then $P \times Q$ is $2^\mu$-c.c). \nl
So a large number are solved.  But, of course, solving two of those 
problems does not necessarily solve their natural combinations.

\newpage

\head {\S1} \endhead  \resetall
\bigskip

We return here to consistency of statements of the form  $\chi \rightarrow  
[\theta]^2_{\sigma,2}$ (i.e. for every  $c:[\chi]^2 \rightarrow \sigma$  
there is  $A \in  [\chi]^\theta$ such that on  $[A],c$ has at most two 
values), (when  $2^\mu \ge \chi > \theta^{< \mu} > \mu$, of course).  
In \cite[\S2]{Sh:276} this was done for $\mu = \aleph_0,\chi = 2^\mu,
\theta = \aleph_1,2 < \sigma < \omega$ and $\chi$ quite large 
(in the original universe $\chi$ is an Erd\"os cardinal).  Originally,
\cite[\S2]{Sh:276} was written
for any $\mu = \mu^{<\mu}$ ($\chi$  measurable in the original universe)
but because of the referee urging it is written up there for $\mu = 
\aleph_0$ only; though with an eye on the more general result which is 
only stated.  In \cite{Sh:288} the main 
objective is to replace colouring of pairs by colouring of $n$-tuples
(and even $(< \omega )$-tuples) but we also say somewhat more on the
$\mu > \aleph_0$ case (in \cite[1.4]{Sh:288}) and using only $k^2_2$-Mahlo 
(for a specific natural number $k^2_2$)(an improvement for $\mu = \aleph_0$ 
too), explaining that it is like \cite{Sh:289}.  
A side benefit of the present paper is giving a full self-contained proof 
of this theorem even for 1-Mahlo.  The main point of this work is to 
increase $\theta$, and this time write it for $\mu = \mu^{< \mu} > \aleph_0$, 
too.
\smallskip

The case $\theta = \mu^+$ is easier as it enables us to separate the forcing
producing the sets admitting few colours: each appear for some $\delta <
\chi,cf(\delta) = \mu^+$, is connected to a closed subset $a_\delta$ of
$\delta$ unbounded in $\delta$ of order type $\mu^+$, so that below 
$\alpha < \delta$ in $P_\alpha$ we
get little information on the colouring on the relevant set.  Here there is
less separation, as names of such colouring can have long common initial
segments, but they behave like a tree and in each node we divide the set to
$\mu$ sets, each admitting only 2 colours. \nl
As we would like to prove the theorem also for $\mu > \aleph_0$, we repeat 
material on $\mu^+$-c.c., essentially from \cite{Sh:80}, \cite{ShSt:154a},
\cite{Sh:288}.
\newpage

\definition{\stag{1.1} Definition}:  1) Let $D$ be a normal filter on
$\mu^+$ to which $\{\delta < \mu^+:\text{cf}(\delta = \mu\}$ belongs.
A forcing notion $Q$ satisfies $*^\epsilon_D$ where $\epsilon$ is a 
limit ordinal $< \mu$, if player I has a winning strategy in the 
following game $*^\epsilon_D[Q]$ defined 
as follows:
\newline
\noindent
\underbar{Playing}: the play finishes after  $\epsilon$ moves.   

In the $\zeta$-th move:   

Player I --- if  $\zeta \ne 0$  he chooses  
$\langle q^\zeta_i:i < \mu^+ \rangle$ such that $q^\zeta_i \in Q$
\roster
\item "{{}}"  $\quad \quad$
and $(\forall \xi < \zeta)(\forall^D i < \mu^+)
p^\xi_i \le q^\zeta_i$ and he chooses a 
\item "{{}}"   $\quad \quad$
function $f_\zeta:\mu^+ \rightarrow \mu^+$ such that for a club of 
$i < \mu^+,f_\zeta(i) < i$;  
\item "{{}}"  $\quad \quad$
if $\zeta = 0$ let  $q^\zeta_i = \emptyset_Q$, $f_\zeta =$ is identically
zero.
\endroster   
\medskip

Player II --- he chooses  $\langle p^\zeta_i:i < \mu^+ \rangle$  
such that $(\forall^D i) q^\zeta_i \le p^\zeta_i$ and $p^\zeta_i \in 
Q$.
\medskip

\noindent
\underbar{The Outcome}:  Player I wins provided that for some $E \in D$: 
if \newline
$\mu < i < j < \mu^+,i,j \in E$,
$cf(i) = cf(j) = \mu$ and $\dsize \bigwedge_{\xi < \epsilon} f_\xi
(i) = f_\xi(j)$ \underbar{then} the set  
$\{p^\zeta_i:\zeta < \epsilon \} \cup \{p^\zeta_j:\zeta <
\epsilon \}$  has an upper bound in  $Q$. \newline
1A) If $D$ is $\{A \subseteq \mu^+:\text{for some club } E \text{ of }
\mu^+ \text{ we have } i \in E \and \text{ cf}(i) = \mu \Rightarrow i \in A\}$
we may write $\mu$ instead of $D$ (in $*^\varepsilon_D$ and in the related
notions defined below and above). \nl
2)  A strategy for a player is a sequence $\bar F = \langle F_\zeta:\zeta < 
\epsilon \rangle,F_\zeta$ telling him what to do in the $\zeta$-th move 
depending only on the previous moves of the other player.  \underbar{But} 
here a play according to the strategy $\bar F$ will mean the player chooses 
in the $\zeta$-th move for each $i < \mu^+$ an element of $Q$ which is 
possibly strictly above (in $\le_Q$'s sense) of what $F_\zeta$ dictates and
a function $f_\zeta$ such that on some $E \in D$, the equivalence relation
$f_\zeta(\alpha) = f_\zeta(\beta)$ induce on $E$ refine the one which 
the strategy induces (this change does not change the truth value 
of ``player $X$ has a winning strategy").  This applies to the game
$\otimes^\varepsilon_Q$ in part (5) below. \newline
3)  We define $**^\varepsilon_\mu$ similarly but for $\zeta$ limit 
$q^\zeta_i$ is not chosen (so player II has to satisfy for limit $\zeta$
just $\forall \xi < \zeta \Rightarrow (\forall^D \,i)(p^\xi_i \le
p^\zeta_i)$).  \nl
4) We may allow the strategy to be non-deterministic, e.g. choose not
$f_\zeta$ just $f_\zeta/D_{\mu^+}$. \nl
5) We say a forcing notion $Q$ is $\varepsilon$-strategically complete if for
the following game, $\bigotimes^\varepsilon_Q$ player I has a winning
strategy. \nl
In the $\zeta$-th move: \nl
Player I - if $\zeta \ne 0$ he chooses $q_\zeta \in Q$ such that
$(\forall \xi < \zeta) p_\xi \le q_\zeta$ if $\zeta = 0$ let $q_\zeta =
\emptyset_Q$. \nl
Player II - he chooses $p_\zeta \in Q$ such that $q_\zeta \le p_\zeta$.
\mn
\ub{The Outcome}:  In the end Player I wins provided that he always has a
legal move.  \nl
6) We say $Q$ is $(< \mu)$-strategically complete if for each
$\varepsilon < \mu$ it is $\varepsilon$-strategically closed. 
\enddefinition
\bigskip

\remark{\stag{1.1A} Remark}  1) In this paper, in the case 
$\mu = \aleph_0$ we can use the
Knaster condition instead of $*^\varepsilon_\mu$. \newline
2) We use below $*^\varepsilon_\mu$ and not $**^\varepsilon_\mu$ but
$**^\varepsilon_\mu$ could serve as well. \nl
3) We may consider omitting the strategic completeness (a weak version of
it is hidden in player $I$ winning $*^\varepsilon_D[Q]$), but no present use.
\endremark
\bigskip

\definition{\stag{1.1B} Definition}  1) Let $\bar F^\ell = 
\langle F^\ell_\zeta:\zeta < \varepsilon \rangle$ be a strategy for player I
in the game $*^\varepsilon_D[Q]$ for $\ell = 1,2$.  We say $\bar F^1 \le
\bar F^2$ equivalently, $\bar F^2$ is above $\bar F^1$ if any play $\langle 
(\bar q^\zeta,f_\zeta,
\bar p^\zeta):\zeta < \varepsilon \rangle$ in which player I uses the
strategy $\bar F^2$ (that is letting $(\langle q'_i:i < \mu^+ \rangle,f) =
F_\zeta(\langle \bar p^\xi:\xi < \zeta \rangle)$ we have $i < \mu^+
\Rightarrow q'_i \le q^\zeta_i$ and for some $E \in D,i \in E \and j \in E
\wedge f(i) = f(j) \Rightarrow f_\zeta(i) = f_\zeta(i)$) is also a play in
which player I uses the strategy $\bar F^1$. \nl
2) Let $\alpha^* < \beta^* < \mu,\bold{St}$ be a winning strategy for player I
in the game $\otimes^\beta_Q$.  We say $\langle \bar F^\alpha:\alpha <
\alpha^* \rangle$ is an increasing sequence of strategies of player I 
in $*^\varepsilon_D [Q]$ obeying $\bold{St}$ if:
\mr
\item "{$(a)$}"  $\bar F^\alpha$ is a winning strategy of player I in
$*^\varepsilon_D[Q]$ 
\sn
\item "{$(b)$}"  for $\alpha < \beta < \alpha^*,\bar F^\beta$ is above
$\bar F^\alpha$
\sn
\item "{$(c)$}"  if $\langle (\bar q^\zeta,f_\zeta,\bar p^\zeta):\zeta <
\varepsilon \rangle$ is a play of $*^\varepsilon_D[Q]$, Player I uses his 
strategy $\bar F^\beta$, \ub{then} for any $i < \mu^+$, letting
$F^\alpha(\langle \bar p^\xi:\xi < \zeta \rangle) = (\bar q^{\alpha,\xi},
f'_{\alpha,\zeta})$ we have:

$$
Q \models \bold{St}(\langle q^{\alpha,\xi}_i:\xi < \zeta \rangle) \le
q^{\alpha,\zeta}_i.
$$
\ermn
3) Similarly to (1), (2) for the game $\otimes^\varepsilon_Q$ (instead
$*^\varepsilon_D[Q]$).
\enddefinition
\bigskip

\demo{\stag{1.1C} Observation}  1) Assume $Q$ is $\mu$-complete.  If $\delta
< \mu$ and $\langle \bar F^\alpha:\alpha < \delta \rangle$ is an increasing
sequence of winning strategies of player I in $*^\varepsilon_D[Q]$,
\ub{then} some winning strategy $\bar F^\delta$ of player I in
$*^\varepsilon_D[Q]$ is above every $\bar F^\alpha (\alpha < \delta)$. \nl
2) Assume $\beta^* < \mu$ and $Q$ is $\beta^*$-strategically complete with
a winning strategy $\bold{St}$.  If $\beta < \beta^*$ and $\langle
\bar F^\alpha:\alpha < \beta \rangle$ is an increasing sequence of winning
strategies of player I in $*^\varepsilon_D[Q]$ obeying $\bold{St}$,
\ub{then} for some $\bar F^\beta,\langle F^\alpha:\alpha < \beta +1 \rangle$
is an increasing sequence of winning strategies of player I in $*^\varepsilon
_D[Q]$ obeying $\bold{St}$. \nl
3) Similarly with $\otimes^\varepsilon_Q$ instead of $*^\varepsilon_Q[D]$.
\enddemo
\bigskip

\demo{Proof}  Straight.
\enddemo
\bigskip

\definition{\stag{1.2} Definition}  Assume  $P,R$  are forcing notions,
$P \subseteq  R$,  $P \lessdot R$. \newline
1) We say $\restriction$ is a restriction operation for the pair $(P,R)$
(or $(P,R,\restriction)$ is a strong restriction triple) \underbar{if} 
($P,Q$ are as above, of course, and) for every member $r \in  R,r 
\restriction P \in P$ is defined such that:
\mr
\item "{(a)}"  $r \restriction P \le r$,
\sn
\item "{(b)}"  if $r \restriction P \le p \in P$ then $r,p$ are compatible 
in $R$ in fact have a lub
\sn
\item "{(c)}"  $r^1 \le r^2 \Rightarrow  r^1 \restriction P \le r^2 
\restriction P$
\sn
\item "{(d)}"  if $p \in P$ then $p \restriction P = p$ and $\emptyset_Q
\restriction P = \emptyset_P$
\ermn

(so this is a strong, explicit way to say $P \lessdot R$). \nl
1A) We say weak restriction triple if we omit in clause (b) the ``have a
lub". \nl
2) We say ``$(P,R,\restriction)$ is $\varepsilon$-strategically complete" if
\mr
\item "{$(\alpha)$}"  $\restriction$ is a restriction operation for the pair
$(P,R)$
\sn
\item "{$(\beta)$}"  $P$ is $\varepsilon$-strategically complete
\sn
\item "{$(\gamma)$}"  if $St_1$ is a winning strategy for player I in the
game, $\bigotimes^\varepsilon_p$, \ub{then} in the game 
$\bigotimes^\varepsilon = \bigotimes^\varepsilon[P,R;St_1]$ 
the first player has a winning strategy $St_2$.
\ermn
\ub{Playing}:  A play of $\bigotimes^\varepsilon$ is a play
$\langle(p_\zeta,q_\zeta):\zeta < \varepsilon \rangle$ of 
$\bigotimes^\varepsilon_R$ \ub{but}
\mr
\item "{$(\alpha)$}"   $\langle (q_\zeta \restriction P,q_\zeta 
\restriction P):\zeta < \varepsilon \rangle$ is a play of the game 
$\bigotimes^\varepsilon_P$ in which the first player uses the strategy 
$St_1$ (see \scite{1.1}(2)!).
\ermn
\ub{Outcome}:  If condition $(\beta)_\zeta$ below fails in stage $\zeta$ 
for some $\zeta < \varepsilon$ then the first player loses immediately, and 
if not, then he wins.
\mr
\item "{$(\beta)_\zeta$}"  for every $\zeta < \varepsilon$, if $p \in P$ is 
above $q_\zeta \restriction P$ then $\{p\} \cup \{q_\xi:\xi < \zeta\}$ has 
an upper bound.  (Read second sentence in \scite{1.1}(2)).
\ermn
2A) We say $(P,R,\restriction)$ is $(< \varepsilon)$-strategically complete
if it is $\zeta$-strategically complete for every $\zeta < \varepsilon$. \nl
3) Let ``$(P,R,\restriction)$ satisfy $\ast^\epsilon_\mu$" mean
(usually $\restriction$ will be understood from context hence
omitted):
\mr
\item "{$(\alpha)$}"  $\restriction$ is a restriction operation for the
pair $(P,R)$
\sn
\item "{$(\beta)$}"   $P$  satisfies  $\ast^\epsilon_\mu$
\sn
\item "{$(\gamma)$}"  If  $St_1$ is a winning strategy for player I in
the game $\ast^\epsilon_\mu[P]$ \underbar{then} in the following game 
called $\ast^\epsilon_\mu[P,R;St_1]$ the first player has a winning 
strategy $St_2$.
\ermn
\underbar{Playing}: As before in $(*)^\varepsilon_\mu[R]$, but  
$\langle < q^\zeta_i \restriction P:i < \mu^+ >,< p^\zeta_i \restriction P:
i < \mu^+ >$, \newline

$\qquad \quad f_\zeta:\zeta < \epsilon \rangle$ is 
required to be a play of $*^\epsilon_\mu[P]$ in which
first player \newline

$\qquad \quad$ uses the strategy $St_1$ 
(see the second sentence of \scite{1.1}(2)).
\medskip
    
\noindent
We also demand that if $\{p^\zeta_j:j < i\} \subseteq P$, then 
$q^\zeta_i \in P$.
\medskip

\noindent
\underbar{The outcome}:  Player I wins provided that:
\mr
\item "{$(*)$}"  for some club $E$ of $\mu^+$ if
$i < j$ are from $E$, cf$(i) = cf(j) = \mu,\dsize 
\bigwedge_{\xi < \varepsilon} f_\xi(i) = f_\xi(j)$  and $r \in P$ is a 
$\le_P$-upper bound of $\{p^\zeta_i \restriction P:\zeta < \epsilon \} 
\cup \{p^\zeta_j \restriction P:\zeta < \epsilon \}$, \underbar{then}
\footnote{could let some strategy determine $r$, no need at present}
$\{r\} \cup \{p^\zeta_i:\zeta < \epsilon \} \cup \{p^\zeta_j:\zeta < 
\epsilon \}$ has an upper bound in $R$.
\ermn
In this case we say that $St_2$ projects to $St_1$ or is above $St_1$.  
If we omit the demand
on the outcome (so maybe $St_2$ is not a winning strategy of player I in
$*^\varepsilon_\mu[R]$), we say $St_2$ weakly projects to $St_1$.
\medskip

\noindent
\underbar{Note}: Naturally in $St_2$ the functions $f_\zeta$ code more
information than $St_1$, we may use a function $g$ to decode the ``older" 
part.  \nl
3A) The game $*^\varepsilon_D[P,R,\restriction]$ and ``$(P,R,\restriction)$
satisfies $*^\varepsilon_D$" are defined naturally and similarly projections
of strategies. \nl
4) We say $(P,R,\restriction)$ satisfies strongly $*^\varepsilon_\mu$ if (when
$\restriction$ is clear from context, it is omitted):
\mr
\item "{$(\alpha)$}"  $\restriction$ is a restriction operation for the
pair $(P,R)$
\sn
\item "{$(\beta)$}"  $P$ satisfies $*^\varepsilon_\mu$
\sn
\item "{$(\gamma)$}"  the first player has a winning strategy in the game
$*^\varepsilon_\mu[P,R,\restriction]$ where
\ermn
\ub{Playing}:  Just like a play of $*^\varepsilon_\mu[R]$, except that
\mr
\item "{$\bigoplus$}"  in addition, for every limit ordinal 
$\zeta < \varepsilon$, in the $\zeta$-the move first the second player
is allowed to choose $\langle r^\zeta_i:i < \mu^+ \rangle$ such that:
$\dsize \bigwedge_{\xi < \zeta} p^\xi_i \restriction P \le r^\zeta_i \in P$
is an upper bound of $\{q^\xi_i \restriction P:\xi < \zeta\}$ and the first
player choosing $q^\zeta_i$ has to satisfy also $(\forall^D i)(r^\zeta_i \le
q^\zeta_i)$.
\ermn
\ub{Outcome}:  Player I wins if $(*)$ from part (3) holds or
\mr
\item "{$(*)^-$}"  in the play $\left< \langle p^\zeta_i \restriction P:
i < \mu^+ \rangle,\langle q^\zeta_i \restriction P:i < \mu^+ \rangle:\zeta
< \varepsilon \right>$ of $*^\varepsilon_\mu[P]$ the first player loses,
(note concerning the outcome, then  now in $(*)$ in part (3), the existence
of $r$ is not (even essentially) guaranteed.) 
\ermn
5) If $\restriction_\ell$ is a 
restriction operation for $(P_\ell,P_{\ell +1})$ for $\ell =1,2,\restriction
= \restriction_1 \circ \restriction_2$, \ub{then} ``a strategy $St$ of first
player in $*^\varepsilon_\mu[P_1,P_3]$ project to one for $*^\varepsilon_\mu
[P_1,P_2]$" is defined naturally.
\enddefinition
\bigskip

\remark{\stag{1.2A} Remark}  We may restrict ourselves to a 
suitable family of strategies $St_1$ (to work in the iteration this 
family has to be suitably closed).
\endremark
\bigskip

\proclaim{\stag{1.3} Claim}  1) If the forcing notion $P$ satisfies
$*^\varepsilon_\mu$ \underbar{then} $P$ satisfies the $\mu^+$-c.c. \newline
2) If $P$ satisfies $*^\varepsilon_\mu$ and $R$ is the trivial forcing
$\{ \emptyset_P \}$ \underbar{then} the pair $(R,P)$ satisfies  
$*^\varepsilon_\mu$ where $\restriction$ is defined by $p \restriction 
R = \emptyset$. \newline
3) If $(P,R,\restriction)$ satisfies $*^\varepsilon_\mu$ \underbar{then} $P$ 
and $R$ satisfy $*^\varepsilon_\mu$. \newline
4) If triples $(P_0,P_1,\restriction_0),(P_1,P_2,\restriction_1)$ satisfy 
$*^\varepsilon_\mu$ \underbar{then} $(P_0,P_2,\restriction_0 \circ 
\restriction_1)$ satisfies $*^\varepsilon_\mu$. \newline
5) If $P$ satisfies $*^\varepsilon_\mu$ and $\Vdash_P$ ``$\underset\tilde
{}\to Q$ satisfies $*^\varepsilon_\mu$" \underbar{then} $P * \underset\tilde
{}\to Q$ satisfies $*^\varepsilon_\mu$ moreover the pair $(P,P *
\underset\tilde {}\to Q)$ (with the natural $\restriction$) satisfies 
$*^\varepsilon_\mu$.
\endproclaim
\bigskip

\demo{Proof}  Should be clear.
\enddemo
\bigskip

\remark{\stag{1.3A} Remark}  1) if $D$ is a normal filter on $\mu^+$ to which
$\{\delta < \mu^+:\text{cf}(\delta) = \mu\}$ belongs, then in \scite{1.3}
we can repalce $*^\varepsilon_\mu$ by $*^\varepsilon_D$ (of course, in
part (5), $D$ in $V^P$ means the normal filter it generates). \nl
Similarly for the claim below. \nl
2) Assume that in the game of choosing $A_i \in D^+$ for $i < \varepsilon$
(or $i < \mu$), with player I choosing $A_{2i}$, player II choosing
$A_{2i+1},A_i$ decreasing, player II loses iff he sometime has no legal move;
player I has a strategy guaranteeing that he has legal moves.
(If $\kappa$ in measurable $V$ in $V^{\text{Levy}(\mu < \kappa)}$ this
holds for some $D$ by \cite{JMMP}.)  In fact assume more generally that
${\Cal P}$ is a partial order and ${\Cal F}:
{\Cal P} \rightarrow \{A:A \subseteq \mu^+\}$ is decreasing: ${\Cal P}
\models x \le y \Rightarrow {\Cal F}(y) \subseteq {\Cal F}(x)$ 
and ${\Cal E}$ is a function with domain ${\Cal P}$ where ${\Cal E}(x)$ 
is a non-empty subset of $[{\Cal F}(x)]^2$
and ${\Cal P} \models x \le y \Rightarrow {\Cal E}(y) \subseteq {\Cal E}(y)$
(above ${\Cal P} = (D^+,\supseteq),{\Cal F}$ is the identity and we say that
a forcing notion $Q$ satisfies
\mr
\item "{$*^\varepsilon_{{\Cal P},{\Cal F},{\Cal E}}$}"  $\quad$ if in 
the following game $*^\varepsilon_{{\Cal P},{\Cal F},{\Cal E}}[Q]$, the 
first player has \nl

$\,\,$ a winning strategy.
\ermn
A play last $\varepsilon$ moves, in the $\zeta$-th move player I chooses
$x_\zeta \in {\Cal P}$ such that $\xi < \zeta \Rightarrow y_\xi \le_{\Cal P}
x_\zeta$ and $\langle q^\zeta_i:i \in {\Cal F}(x_\zeta) \rangle$ such that
$\xi < \zeta \and i \in {\Cal F}(x_\zeta) \Rightarrow p^\xi_i \le q^\zeta_i$
and player II chooses $y_\zeta \in {\Cal P}$ such that $x_\zeta \le
y_\zeta$ and $\langle p_i:i \in {\Cal F}(y_\zeta) \rangle$ such that $i \in
{\Cal F}(y_\zeta) \Rightarrow q^\zeta_i \le_Q P^\zeta_i$.
\sn
\ub{Outcome}:  Player I wins a play if
\mr
\item "{$(\alpha)$}"  for every limit $\zeta < \varepsilon$ he has a legal
move (this depends on having upper bounds in ${\Cal P}$ and in $Q$)
\sn
\item "{$(\beta)$}"  for every $\{i,j\} \in \dbca_{\zeta < \varepsilon}
{\Cal E}(x_\zeta)$, in $Q$ there is an upper bound to \nl
$\{p^\zeta_i:\zeta < \varepsilon\} \cup \{p^\zeta_j:\zeta < \varepsilon\}$.
\ermn
The natural generalizations of the relevant lemmas works for this notion. \nl
3)  We can systematically use the weak restriction triples, and/or use the
strong version of $*^\varepsilon_\mu$ for triples in this paper.
\endremark
\bigskip

\proclaim{\stag{1.4} Claim}  1) If the forcing notions $P_1,P_2$ are 
equivalent \underbar{then} $P_1$ satisfies $*^\varepsilon_\mu$ iff 
$P_2$ satisfies $*^\varepsilon_\mu$. \newline
2) Suppose $\restriction$ is a restriction operation for $(P_1,P_2),B_\ell$
the complete Boolean algebra corresponding to $P_\ell$ (so $B_1 \lessdot
B_2)$ and $\restriction'$ is the projection from $B_2$ to $B_1$ and
$P'_\ell = (B_\ell \backslash \{0\},\ge)$ then
\medskip
\roster
\item "{$(a)$}"  $(P'_1,P'_2,\restriction')$ is a restriction triple and
\sn
\item "{$(b)$}"  $(P_1,P_2,\restriction)$ satisfies $*^\varepsilon_\mu$ 
iff $(P'_1,P'_2,\restriction')$ satisfies $*^\varepsilon_\mu$.
\endroster
\smallskip
\noindent
2A)  In part (2) it is enough to assume that $\restriction$ is a weak
restriction operation. \nl
3) If a forcing notion $Q$ satisfies 
$*^\varepsilon_\mu$ \underbar{then} player I has a
winning strategy in the play even if we demand from him: $\dsize \bigwedge
_{\xi < \zeta} p^\xi_i = \emptyset_Q \Rightarrow q^\zeta_i = \emptyset_Q$ 
for each $i < \mu^+$. \newline
4) Similarly for $(P,R,\restriction)$ satisfying $*^\varepsilon_\mu$ 
demanding $\dsize \bigwedge_{\xi < \zeta} p^\xi_i = \emptyset_R \Rightarrow
q^\zeta_i = \emptyset_R$ and $\dsize \bigwedge_{\xi < \zeta} p^\xi_i 
\Rightarrow q^\zeta_i \in P$.
\endproclaim
\bigskip

\demo{\stag{1.4A} Convention}  Strategies are as in \scite{1.4}(3),(4).
\enddemo
\bigskip

\proclaim{\stag{1.5} Definition/Claim}  Assume for $\ell = 1,2$ that 
$(P,R_\ell,\restriction_\ell)$ is a restriction triple, 
$(P,R_\ell,\restriction_\ell)$ satisfies $*^\varepsilon_\mu$,
and we let \newline
$R = \{(p,r_1,r_2):p \in P,r_1 \in R_1,r_2 \in R_2,P \models
``r_1 \restriction P \le p"$ and $P \models ``r_2 \restriction P \le p"\}$
\newline
identifying $r_1 \in R_1$ with $(r_1 \restriction P,r_1,\emptyset_{R_2})$,
and identifying $r_2 \in R_2$ with \newline
$(r_2 \restriction P,\emptyset_{R_1},r_2)$. \nl
Under the quasi order

$$
\align
(p,r_1,r_2) &\le (p'_1,r'_1,r'_2) \text{ iff } p \le_P p' \\
  &\and \text{ lub}_{R_1}\{p,r_1\} \le_{R_1} \text{ lub}_{R_1}\{p,r'_1\} \\
  &\and \text{ lub}_{R_2}\{p,r_2\} \le_{R_2} \text{ lub}_{R_2}\{r''_2\}.
\endalign
$$
\mn
\underbar{Then} $R_\ell \lessdot R$ (for $\ell = 1,2$) and
$(R_\ell,R,\restriction'_\ell)$ is a restriction triple and it
satisfies $*^\varepsilon_\mu$, where 
$(p,r_1,r_2) \restriction^\prime_\ell R_\ell =$ the lub of $p,r_\ell$ in
$R_\ell$ (see clause (b) of Definition \scite{1.2}(1)).
\endproclaim
\bigskip

\definition{\stag{1.6} Definition/Lemma}  Let  $\mu = \mu^{<\mu} < \kappa = 
\text{ cf}(\kappa) \le \lambda \le \chi$.  (Usually fixed hence suppressed 
in the notation).  We define and prove the following by induction on (the 
ordinal) $\alpha$: \newline
1) [Def].  Let  ${\Cal K}^\alpha = {\Cal K}^\alpha_{\mu,\kappa,\lambda,
\chi}$ be the family of 
sequences  $\bar Q = \langle P_\beta,{\underset\tilde {}\to Q_\beta},a_\beta:
\beta < \alpha \rangle$  such that:
\medskip
\roster
\item "{(a)}"  $\langle P_\beta,{\underset\tilde {}\to Q_\beta}:\beta
< \alpha \rangle$  is a $(< \mu)$-support iteration (so $P_\alpha = 
\text{ Lim}_\mu \bar Q$  denotes the natural limit)
\sn
\item "{(b)}"  $a_\beta \subseteq \beta$, $|a_\beta| < \kappa$,
$[\gamma \in a_\beta \Rightarrow  a_\gamma \subseteq a_\beta]$
\sn
\item "{(c)}" ${\underset\tilde {}\to Q_\beta}$ is strategically 
$(< \mu)$-complete, has cardinality $< \lambda$ and is a $P^*_{a_\beta}$-name
(see parts \scite{1.6}(2)(b) and \scite{1.6}(5)(b) below).
\ermn
1A) [Def]  $\bar Q$ is called standard \underbar{if}: for every $\beta <
\ell g(\bar Q)$ each element of
${\underset\tilde {}\to Q_\beta}$ is from $V$, even from ${\Cal H}(\chi)$, 
and the order is a fixed quasi order from $V$ such that any chain of length 
$< \mu$ which has an upper bound has a lub (we can use less), but note that
the set of elements is not necessarily from $V$.

\medskip

\noindent
2) [Def].  For  $\bar Q$  as above:
\medskip
\roster
\item "{(a)}" $a \subseteq \alpha$  is called $\bar Q$-closed if
$[\beta \in  a \Rightarrow  a_\beta \subseteq a]$; we also call it \newline
$\langle a_\beta:\beta < \alpha \rangle$-closed and let $\bar a^{\bar Q} =
\langle a_\beta:\beta < \alpha \rangle$
\smallskip
\noindent
\item "{(b)}"  for a $\bar Q$-closed subset $a$ of $\alpha$  we let
\endroster

$$
\align
P_a = \{ p \in P_\alpha:&\text{Dom}(p) \subseteq a \text{ and for each }
\beta \in \text{ Dom}(p) \\
  &\text{we have: } p(\beta) \text{ is a } P_{a \cap \beta} \text{-name} \\
  &\text{(i.e. involving only } G_{P_\beta} \cap P_{a \cap \beta} \\
  &\text{so necessarily } Q \in V[G_{P_\beta} \cap P_{a \cap \beta}])\}
\endalign
$$

$$
\align
P^*_a = \{ p \in P_\alpha:&\text{Dom}(p) \subseteq a \text{ and for each }
\beta \in \text{ Dom}(p) \text{ we have: } p(\beta) \\
  &\text{ is a } P^*_{a_\beta} \text{-name} \text{ and: if }
{\underset\tilde {}\to Q_\beta} \subseteq V \text{ and } \bar Q
\text{ is standard, then} \\
  &\,p(\beta) \text{ is from } V\}.
\endalign
$$ 
\noindent
On both $P_a$ and $P^*_a$, the order is inherited from $P_\alpha$. 
Note that $P^*_a$ is defined by induction on $\sup(a)$.
\newline
\medskip

\noindent
3) [Lemma]  For  $\bar Q$  as above,  $\beta < \alpha$
\medskip
\roster
\item "{(a)}"  $\bar Q \restriction \beta \in {\Cal K}^\beta$
\sn
\item "{(b)}"  if $a \subseteq \beta$ then: $a$ is $\bar Q$-closed
\underbar{iff} $a$ is $(\bar Q \restriction \beta)$-closed
\sn
\item "{(c)}"  if $a \subseteq \alpha$ is $\bar Q$-closed, \ub{then} so is 
$a \cap \beta$, in fact $\beta$ is $\bar Q$-closed and the intersection of
a family of $\bar Q$-closed subsets of $\alpha$ is $\bar Q$-closed.  
\ermn
4) [Lemma].  For $\bar Q$ as above, and $\beta < \alpha$,
\mr
\item "{$(a)$}"  $P_\beta \lessdot P_\alpha$, moreover, if $p \in P_\alpha,
p \restriction \beta \le q \in P_\beta$ then $(p \restriction 
(\alpha \backslash \beta)) \cup q \in  P_\alpha$ is a lub of $p,q$
\sn
\item "{$(b)$}"  $P_\alpha/P_\beta$ is strategically $(< \mu)$-complete (hence
does not add new sequences of length $< \mu$ of old elements).
\ermn
5) [Lemma].  For  $\bar Q$  as above
\medskip
\roster
\item "{(a)}"  $P^\ast_\alpha$ is a dense subset of  $P_\alpha$
\sn
\item "{(b)}"  if $a$  is $\bar Q$-closed \underbar{then} $P_a \lessdot
P_\alpha$ and $P^\ast_a$ is a dense subset of  $P_a$.
\sn
\item "{(c)}"  if $a$ is $\bar Q$-closed, $p \in P_\alpha,p \restriction a
\le q \in P_a$ \underbar{then} $(p \restriction (\alpha \backslash a)) 
\cup q$ belongs to $P_\alpha$ and is a lub of $p,q$ in $P_\alpha$
\sn
\item "{(d)}"  if $a$ is $\bar Q$-closed, \underbar{then} 
$\bar Q \restriction a \in {\Cal K}^{\text{otp}(a)}$ (up 
to renaming of indexes)
\sn
\item "{$(e)$}"  if $a \subseteq b \subseteq \ell g(\bar Q)$ are 
$\bar Q$-closed, \ub{then} $(P^*_a,P^*_b,\restriction)$ is a restriction 
triple (where $p \restriction P^*_b = p \restriction a$)
\endroster
\medskip 

\noindent
6) [Lemma].  The sequence
$\bar Q = \langle P_\beta,{\underset\tilde {}\to Q_\beta},a_\beta:
\beta < \alpha \rangle$  belongs to  ${\Cal K}^\alpha$ if
$\alpha$  is a limit ordinal and $\dsize \bigwedge_{\gamma < \alpha} \bar Q
\restriction \gamma \in {\Cal K}^\gamma$.
\medskip

\noindent
7) [Lemma].  The sequence  $\bar Q = \langle P_\beta,
{\underset\tilde {}\to Q_\beta},a_\beta:\beta < \alpha \rangle$
belongs to ${\Cal K}^\alpha$ if $\alpha = \gamma + 1$, \newline
$a_\gamma \subseteq \gamma$ is a $(\bar Q \restriction \gamma)$-closed 
set of cardinality $< \kappa,{\underset\tilde {}\to Q_\gamma}$ is a 
$P^*_{a_\gamma}$-name of a \newline
$(< \mu)$-strategically complete forcing notion of cardinality $< \lambda$.
\medskip

\noindent
8) [Def].  ${\Cal K}^{< \alpha} = \dsize \bigcup_{\beta < \alpha} {\Cal K}
^\beta$.
\enddefinition
\bigskip

\demo{Proof}  Straightforward.
\enddemo
\bigskip

\definition{\stag{1.7} Definition}  Let $\mu  = \mu^{< \mu} < \kappa = 
\text{ cf}(\kappa)\le \lambda \le \chi$ (usually fixed hence suppressed 
in the notation) and $\varepsilon$ a limit ordinal $< \mu$.  We define 
the following by induction
on (the ordinal) $\alpha$: \newline
1)  We let ${\Cal K}^{\varepsilon,\alpha} =
{\Cal K}^{\varepsilon,\alpha}_{\mu,\kappa,\lambda,\chi}$ be the family of
sequences \newline
$\bar Q = \langle P_\beta,{\underset\tilde {}\to Q_\beta},a_\beta,
I_\beta:\beta < \alpha \rangle$ such that:
\medskip
\roster
\item "{$(\alpha)$}"  $\langle P_\beta,{\underset\tilde {}\to Q_\beta},
a_\beta:\beta < \alpha \rangle \in {\Cal K}^\alpha$
\sn
\item "{$(\beta)$}"  $I_\beta$ is a family of $\bar Q$-closed (see part (2)
below, it is not what was defined in \scite{1.6}(2)(a))
subsets of $a_\beta$, closed under finite unions, increasing unions of 
length $< \mu$ and such that $\emptyset \in I_\beta$
\sn
\item "{$(\gamma)$}"  each $a_\beta$ is $(\bar Q \restriction \beta)$-closed
(see part (2) below, this is not as in \scite{1.6})
\sn
\item "{$(\delta)$}"  if $b \in I_\beta$ then the pair 
$(P^*_b,P^*_{a_\beta \cup \{ \beta \}})$ satisfies $*^\varepsilon_\mu$,
of course, for the natural restriction operation.
\endroster
\medskip

\noindent
(2) For $\bar Q \in {\Cal K}^{\varepsilon,\alpha}$ (even satisfying just 
\scite{1.7}(1)($\alpha$) + ($\beta$)) we
say that a set $a$ is $\bar Q$-closed in $b$ (or is 
$\langle a_\beta,I_\beta:\beta < \alpha \rangle$-closed) if $a \subseteq b 
\subseteq \alpha,[\beta \in a \Rightarrow
a_\beta \subseteq a]$ and $[\beta \in b \backslash a \Rightarrow a 
\cap a_\beta \in I_\beta]$.  If we omit ``in $b$" we mean $b = \alpha$.
\medskip

\noindent
(3) 
\roster
\item "{(a)}"  $\bar Q$ is simple if for all $\beta < \alpha$ \newline
$I_\beta = \{ b \subseteq a_\beta:b \text{ is } \bar a^{\bar Q}$-closed and: 
for every $\gamma \in a_\beta \cup \{ \beta\}, \text{ if  cf}(\gamma) = \mu^+$
\newline
and $\gamma = \text{ sup}(\gamma \cap b), \text{ then } \gamma \in b\}$.
\smallskip
\noindent
\item "{(b)}"  $\bar Q^- = \langle P_\beta,{\underset\tilde {}\to Q_\beta},
a_\beta:\beta < \alpha \rangle,a^{\bar Q} = \langle a_\beta:\beta < \alpha
\rangle$, and \newline
$\bar I^Q = \langle I_\beta:\beta < \alpha \rangle$  
\smallskip
\noindent
\item "{(c)}"  $\bar Q$ is standard if $\bar Q^-$ is standard
\smallskip
\noindent
\item "{(d)}"  ${\Cal K}^{\varepsilon,< \alpha} = \dsize \bigcup
_{\beta < \alpha} {\Cal K}^{\varepsilon,\beta}$.
\endroster
\enddefinition
\bigskip

\proclaim{\stag{1.8} Claim}  Let $\bar Q \in 
{\Cal K}^{\varepsilon,\alpha}$. \newline
1)  If $\beta < \alpha$ then $\bar Q \restriction \beta =: \langle P_\gamma,
{\underset\tilde {}\to Q_\gamma},a_\gamma,I_\gamma:\gamma < \beta \rangle$ 
belongs to ${\Cal K}^{\varepsilon,\beta}$; moreover, if \newline
$b \subseteq \alpha$
is $\bar a^{\bar Q}$-closed \underbar{then} $\bar Q \restriction b \in 
{\Cal K}^{\varepsilon,\text{otp}(b)}$ (up to renaming of index sets)
understanding $I^{\bar Q \restriction b}_\beta = I^{\bar Q}_\beta \restriction
b$. \newline
2)  If $a \subseteq b \subseteq \beta \le \alpha$ and $a$ is $\bar Q$-closed 
in $b$ \ub{then}:  $a$ is ($\bar Q \restriction \beta$)-closed in $b$. \nl
3)  If $\beta < \alpha,a \subseteq \alpha$ is $\bar Q$-closed and $\gamma
\in \alpha \backslash \beta \Rightarrow a \cap a_\gamma \in I_\gamma$,
\ub{then} $a \cap \beta$ is $\bar Q$-closed. \newline
4) If $\bar Q$ is simple, $\beta < \alpha,a \subseteq \alpha$ is 
$\bar Q$-closed and cf$(\beta) \ne \mu^+ \vee (\forall \gamma \in \alpha
\backslash \beta)(a_\gamma \cap a \cap \beta$ is bounded in $\beta$),
\ub{then} $a \cap \beta$ is $\bar Q$-closed. \nl
5)  The family of $\bar Q$-closed $a \subseteq \alpha$ is closed under
increasing union of length $< \mu$ and $\emptyset$ belongs to it. \newline
6)  If $a,b$ are $\bar Q$-closed, then so is $a \cup b$. \nl
7)  If $a \subseteq b \subseteq c \subseteq \ell g(\bar Q),a$ is 
$\bar Q$-closed in $c$,
\ub{then} $a$ is $\bar Q$-closed in $b$. \nl
8) If $a \subseteq b \subseteq \alpha,a$ is $\bar Q$-closed in $b$,
\ub{then} $a \cap \alpha$ is $(Q \restriction \beta)$-closed in 
$b \cap \beta$.
\endproclaim
\bigskip

\demo{Proof}  Straight.
\enddemo
\bigskip

\remark{\stag{1.8A} Remark}  Simple $\bar Q$ is what we shall use.
\endremark
\bigskip

\proclaim{\stag{1.9} Lemma}  Assume $\bar Q \in {\Cal K}^{\varepsilon,\alpha}$
and $a,b$ are $\bar Q^-$-closed subsets of $\alpha$ and $a$
is a $\bar Q$-closed subset of $b \, (\subseteq \alpha)$ and $\bar Q$ is
simple or at least
\medskip
\roster
\item "{$(*)$}"  $\gamma < \beta < \alpha \Rightarrow a_\beta \cap
(\gamma + 1) \in I_\beta$ \newline
(hence $\gamma < \beta < \alpha \and \text{cf}(\gamma) < \mu \Rightarrow
a_\beta \cap \gamma \in I_\beta)$.
\endroster
\medskip
\noindent
\underbar{Then} the pair $(P^*_a,P^*_b)$ satisfies 
$\ast^\varepsilon_\mu$.
\endproclaim
\bigskip

\demo{Proof}  We can assume by \scite{1.8}(1) that $b = \alpha$.  
By induction on $\alpha$ we shall show that for all $\bar Q$-closed 
subsets $a$ of $\alpha$ the pair $(P^*_a,P^*_\alpha)$ satisfies 
$*^\varepsilon_\mu$ (see Definition \scite{1.2}(3)) and this is proved first
when $a = \emptyset$ and then when $a \ne \emptyset$.  
So we fix a strategy $St_a$ for the first
player in $*^\varepsilon_\mu[P^*_a]$; why it exists?  If $a = \emptyset$,
trivially, if $a \ne \emptyset$ by the way the proof is arranged we know
the conclusion for $(a',b') = (\emptyset,a)$, and as otp$(a) \le \alpha$
clearly $St_a$ exists.  Next we shall choose a strategy for 
the first player in the game $*^\varepsilon_\mu[P^*_a,P^*_\alpha,St_a]$, 
where at stage $\zeta < \varepsilon$ the first player chooses 
$\{ q^\zeta_\xi:\xi < \mu^+\}$, a regressive function $f_\zeta$ 
from $\mu^+$ to $\mu^+$ and the second player replies with suitable 
$\{ p^\zeta_\xi:\xi < \mu^+\}$.

For simplicity the reader may assume that the 
${\underset\tilde {}\to Q_\beta}$ are
$\mu$-complete (which is the case used; otherwise we have to use the
$(< \mu)$-strategic completeness (and remember \scite{1.1}(2) second
sentence).
\bigskip

\noindent
\underbar{Case 1}:  $\alpha = \beta + 1,\beta \in a$.

So $a_\beta \subseteq a$, now $a \cap \beta$ is $(\bar Q \restriction 
\beta)$-closed (by \scite{1.8}(2)) hence by the induction hypothesis 
$(P^*_{a \cap \beta},P^*_\beta)$ satisfies $*^\varepsilon_\mu$. 
Apply \scite{1.5} with $P^*_{a \cap \beta},P^*_\beta,P^*_a$ here standing for
$P,R_1,R_2$ there and we get that $(R_2,R)$ satisfies $*^\varepsilon_\mu$,
which (translating) is the desired conclusion.
\bigskip

\noindent
\underbar{Case 2}:  $\alpha = \beta + 1, \beta \notin a$.

We know that $a \cap a_\beta \in I_\beta$. \newline
By Definition \scite{1.7}(1)$(\delta)$ we know that $(P^*_{a \cap a_\beta},
P^*_{a_\beta \cup \{ \beta \}})$ satisfies $*^\varepsilon_\mu$.  By 
\scite{1.5} we get that $(P^*_a,P^*_{a_\beta \cup \{ \beta \} \cup a})$ 
satisfies $*^\varepsilon_\mu$.  Now $a' =: a_\beta \cup \{ \beta \} \cup a$ is
$\bar Q$-closed by \scite{1.8}(6) and $\beta \in a'$ so by Case 1 we have:
$(P^*_{a'},P^*_\alpha)$ satisfies $*^\varepsilon_\mu$.  
Together by \scite{1.3}(4) we have: $(P^*_a,P^*_\alpha)$ satisfies 
$*^\varepsilon_\mu$.
\bigskip

\noindent
\underbar{Case 3}:  $\alpha$ a limit ordinal, cf$(\alpha) \le \mu$.

Here we use \scite{1.4}(3) (i.e. \scite{1.4}(A)). \newline
We can find an increasing continuous sequence $\langle \gamma_\Upsilon:
\Upsilon < \text{ cf}(\alpha) \rangle$ of ordinals $< \alpha$ with limit
$\alpha,\gamma_0 = 0$ and $\gamma_{\Upsilon + 1}$ a successor ordinal.  
Note that
$(a \cap \gamma_{\Upsilon + 1}) \cup \gamma_\Upsilon$ is \newline
$(\bar Q \restriction
\gamma_{\Upsilon + 1})$-closed as $[\gamma_\Upsilon \text{ limit }
\Rightarrow \Upsilon \text{ limit } \and \text{cf}(\Upsilon) < \mu]$ moreover
$a \cup \Upsilon_\gamma$ is $\bar Q$-closed.  
We define by induction on $\Upsilon \le \text{ cf}(\alpha)$ a strategy 
$St^*_\Upsilon$ of player I in the game
$*^\varepsilon_\mu[P^*_a,P^*_{a \cup \gamma_{\Upsilon}}]$ such that for
$\Upsilon_1 < \Upsilon$ we have that $St^*_\Upsilon$ projects to 
$St^*_{\Upsilon_1}$ (see Definition \scite{1.2}(4)) and $St^*_0$ is 
$St_a$. \nl
If we do not assume that all the ${\underset\tilde {}\to Q_\beta}$ are
$\mu$-complete, then we demand that, moreover, they satisfy:
\mr
\item "{$\boxtimes$}"  if $\left< \langle q^\zeta_i:i < \mu^+ \rangle,
f_\zeta,\langle p^\zeta_i:i < \mu^+ \rangle:\zeta < \varepsilon \right>$ is
a play of $*^\varepsilon_\mu[P^*_a,P^*_{a \cup \gamma \Upsilon},St_a]$,
\ub{then} for any ordinal $\beta$, looking at $\langle q^\zeta_i(\beta),
p^\zeta_i(\beta):\zeta < \varepsilon \rangle$ letting $\zeta(\beta,
\emptyset) = \text{ Min}\{\zeta:q^\zeta_i*\beta) \ne \emptyset_Q\}$ if
$\zeta \in [\zeta(\beta,0),\zeta(\beta,1))$ and $q^\zeta_i \restriction
\beta$ forces that $\langle q^\xi_i(\beta):\xi \in [\zeta,(\beta,0),\zeta] 
\rangle$ is increasing, then $q^\zeta_i \restriction \beta$ forces that some
$\langle q'_\xi,p'_\xi:\xi < \zeta - \zeta(\beta,0) + 1 \rangle$ is a play
of $\otimes^\varepsilon_{\underset\tilde {}\to Q_\beta}$ in which player I
uses a fix winning strategy (as in \scite{1.1}(2)!) and $p'_0 =
q^{\zeta(\beta,0)}_i(\beta)$, (remember $q'_0$ not chosen) and
$0 < \xi < \zeta - \zeta(\beta,0)+1 \Rightarrow q'_\xi = 
q^{\zeta(\beta,0) +\xi}_i(\beta)$ and $0 < \xi < \zeta - \zeta(\beta,0)
\Rightarrow p'_\xi = p^\xi_i(\beta)$.
\ermn
This, of course, puts on us a burden also in successor $\gamma$ just to
increase the condition. \nl
The inductive step is done by \scite{1.5}, the limit stage is straight
(using $\boxtimes$ to show we can).
\bigskip

\noindent
\underbar{Case 4}:  $\alpha$ limit ordinal, cf$(\alpha) > \mu^+$.

During the play, player I in the $\zeta$-th move also chooses an ordinal
$\gamma_\zeta,\gamma_\zeta$ increases continuously with $\zeta,\gamma_0 = 0$
as follows:

$$
\gamma_{\zeta + 1} = \text{ min}\{ \gamma < \alpha:(\forall i < \mu^+)
(\forall \xi \le \zeta)(p^\xi_i,q^\xi_i \in P_\gamma)\}
$$
\medskip

\noindent
and he will make $q^\zeta_i \in P_{\gamma_\zeta}$, and the rest is as in
Case 3. 
\bigskip

\noindent
\underbar{Case 5}:  cf$(\alpha) = \mu^+$.

Let $\langle \gamma_\Upsilon:\Upsilon < \mu^+ \rangle$ be increasing
continuously with limit $\alpha, \gamma_0 = 0$, cf$(\gamma_\Upsilon) \le
\mu$, and we imitate Case 4, separating to different plays according to the
value of \newline
$j^\zeta_i = \text{Min}\{j < i:\text{for each } \xi < \zeta
\text{ we have } p^\xi_i \restriction \gamma_i \in P_{\gamma_j} \text{ and }
q^\xi_i \restriction \gamma_i \in P_{\gamma_j}\}$.
\hfill$\square_{\scite{1.9}}$
\bigskip
\enddemo

\proclaim{\stag{1.10} Claim}  Assume
\medskip
\roster
\item "{(a)}"  $\bar Q = \langle P_\alpha,{\underset\tilde {}\to Q_\alpha},
a_\alpha,I_\alpha:\alpha < \delta \rangle$
\item "{(b)}"  $\delta$ a limit ordinal
\item "{(c)}"  for every $\alpha < \delta$ we have $\bar Q \restriction
\alpha \in {\Cal K}^{\varepsilon,\alpha}$.
\endroster
\medskip

\noindent
\underbar{Then}  $\bar Q \in {\Cal K}^{\varepsilon,\delta}$.
\endproclaim
\bigskip

\demo{Proof}  Check.
\enddemo
\bigskip

\proclaim{\stag{1.11} Claim}  Assume
\medskip
\roster
\item "{(a)}"  $\bar Q \in {\Cal K}^{\varepsilon,\alpha}$
\sn
\item "{(b)}"  $a_\alpha \subseteq \alpha$ is $\bar Q$-closed,
$|a_\alpha| < \kappa$
\sn
\item "{(c)}"  $I_\alpha \subseteq \{ b \subseteq a_\alpha:b \text{ is }
\bar Q$-closed \}
\sn
\item "{(d)}"  $I_\alpha$ is closed under finite unions, $I_\alpha$ is closed
under increasing unions of length $< \mu$ and $\emptyset \in I_\alpha$
\sn
\item "{(e)}"  ${\underset\tilde {}\to Q_\alpha}$ is a
$P^*_{a_\alpha}$-name of a forcing notion of cardinality $< \lambda$
\sn
\item "{(f)}"  if $b \in I_\alpha$ then $(P_b,P^*_{a_\alpha} *
{\underset\tilde {}\to Q_\alpha})$ satisfies $\ast^\varepsilon_\mu$
\sn
\item "{(g)}"  $P_\alpha = \text{ Lim}_\mu \bar Q$.
\endroster
\medskip

\noindent
\underbar{Then} $\bar Q \char 94
\langle P_\alpha,{\underset\tilde {}\to Q_\alpha},a_\alpha,I_\alpha \rangle$
belongs to ${\Cal K}^{\varepsilon,\alpha + 1}$.
\endproclaim
\bigskip

\demo{Proof}  Check.
\enddemo
\bigskip

\proclaim{\stag{1.12} Theorem}  Suppose $\mu = \mu^{<\mu} < \kappa = 
\lambda < \chi$ and $\chi$ is measurable. \nl
1)  For some forcing notion $P$ of cardinality $\chi$,
$\mu$-complete not collapsing cardinalities not changing cofinalities we have:
\newline
$\Vdash_P$ ``$2^\mu = \chi$ and for every $\sigma < \mu$  and
$\theta < \kappa$ we have $\chi \rightarrow [\theta]^2_{\sigma,2}$"
(and for a fixed $\varepsilon$ the Axiom: if $Q$ is a $\mu$-complete
forcing notion of cardinality $< \kappa$ satisfying $\ast^\varepsilon_\mu$ and
${\Cal I}_\alpha \subseteq Q$ dense for $\alpha < \alpha^* < \kappa$
\underbar{then} some directed $G \subseteq Q$ is not disjoint to any
${\Cal I}_\alpha$). \nl
2) We can replace ``$\mu$-complete" by ``$(< \mu)$-strategically complete"
(in the demand on $P$ and, in the axiom, on $Q$. 
\endproclaim
\bigskip

\remark{\stag{1.12A} Remark}  We can add ``$P$ satisfies $*^\varepsilon_\mu$"
if the appropriate squared diamond holds which is true in reasonable inner 
models.
\endremark
\bigskip

\demo{Proof}  We concentrate on part (2).  If
we would like to do part (1), we should just demand all the
${\underset\tilde {}\to Q_i}$ are $\mu$-complete.
\mn
\ub{Stage A}:  Fix $\varepsilon  < \mu$ and let ${\Cal K}^\alpha_* = 
\{ \bar Q \in {\Cal K}^{\varepsilon,\alpha}:\bar Q$ is simple and 
standard$\}$,
${\Cal K}_* = \dsize \bigcup_{\alpha < \chi} {\Cal K}^\alpha_*$.  (Note:
$\bar Q$-closed will mean as in \scite{1.7}(3)(a),\scite{1.7}(2).)
By preliminary forcing without loss of generality  
``$\chi$  measurable" is preserved by forcing with $({}^{\chi >}2,
\trianglelefteq)$ (= adding a Cohen subset of $\chi$), see Laver \cite{L}.  
Let us define a forcing notion  $R$: \newline
$R = \{\bar Q:\bar Q \in {\Cal K}^\alpha_*$ for some $\alpha < \chi$ and 
$\bar Q \in {\Cal H}(\chi) \}$ \newline
ordered by:  $\bar Q^1 \le \bar Q^2$ iff  $\bar Q^1 = \bar Q^2 \restriction
\ell g(\bar Q^1)$. \newline
As  $R$  is equivalent to  $({}^{\chi >}2,\trianglelefteq)$  we know that
in  $V^R$,  $\chi$  is still measurable.  Let  $\bar Q^\chi = \langle
P_\beta,{\underset\tilde {}\to Q_\beta},a_\beta:\beta < \chi
\rangle$ be $\bigcup G_R$ and  $P_\chi$ be the limit so $P^* = P^\ast_\chi 
\subseteq P_\chi$ is a dense subset, those are $R$-names.  Now $R * 
{\underset\tilde {}\to P^\ast}$ is the forcing  $P$  we have promised.
The non-obvious point is  
$\Vdash_{R * {\underset\tilde {}\to P^\ast_\chi}}$ 
``$\chi \rightarrow [\theta]^2_{\sigma,2}$" (where  
$\theta < \kappa,\sigma < \mu$).  So suppose
$(r^*,{\underset\tilde {}\to p^*}) \in R \ast 
{\underset\tilde {}\to P^\ast_\chi}$ and
$(r^*,{\underset\tilde {}\to p^*}) \Vdash$ ``the colouring 
${\underset\tilde {}\to \tau}:[\chi]^2
\rightarrow \sigma$  is a counterexample".
Let  $\chi_1 = (2^\chi)^+$.  Let  $G_R \subseteq  R$  be generic over  $V$,  
$r^\ast  \in  G_R$.  By \cite{Sh:289}, but the meaning is explained below 
in $V^R$ we can find an end extension
strong $(\chi_1,\chi,\chi,2^{\kappa + \lambda + 2^\mu},(\kappa + \lambda
+ 2^\mu)^+,\omega)$-system  $\bar M = \langle M_s:s \in [B]^{<\aleph_0}
\rangle$ such that $M_s \prec ({\Cal H}(\chi_1)^{V[G_R]},{\Cal H}(\chi_1),
\in)$, for $x = \{ \chi,G_R,p^*,{\underset\tilde {}\to \tau} \}$, (i.e. 
$x \in \dsize \bigcap_s M_s$ and $B \in [\chi]^\chi$).  We do not define this
as for helping to prove the next theorem (1.13) we assume less in
$V[G_R],M_s \prec ({\Cal H}(\chi_1)^{V[G_R]},\in,{\Cal H}(\chi_1),G_R)$ and:
\medskip
\roster
\item "{$(*)_0$}"  $\bar M = \langle M_s:s \in [B]^{<(1+n^*)} \rangle$ 
is an end extension $(\chi_1,\chi,\chi,2^{\kappa + \lambda + 2^\mu}$, \nl
$(\kappa + \lambda + 2^\mu)^+,n^*)$-system for
$x$, for some $2 \le n^* \le \omega$.
\ermn
where $(*)_0$ means:
\mr
\item "{$(*)'$}"   $B \in [\chi]^\chi$ and $M_s \prec ({\Cal H}(\chi_1),
\in),x \in \dsize \bigcap_s M_s,M_s \cap M_t = M_{s \cap t}$.  \nl
Furthermore, $\|M_s\| = 2^{\kappa + \lambda+ 2^\mu}$
and $[M_s]^{\kappa + \lambda + 2^\mu} \subseteq M_s$.  In addition, for
$v_1,v_2 \in [B]^n,n < 1 + n^*$ there is $f_{v_1,v_2}$, the unique 
isomorphism from $M_{v_1}$ onto $M_{v_2}$, and:
$|v_1 \cap \varepsilon_1| = |v_2 \cap \varepsilon_2|,\varepsilon_1 \in v_1,
\varepsilon_2 \in v_2 \Rightarrow f_{v_1,v_2}(\varepsilon_1) =
\varepsilon_2$.  Finally, $s \triangleleft t \Rightarrow M_s \cap \chi
\triangleleft M_t \cap \chi$.
\ermn
We meanwhile concentrate on case $n^* = 2$.
\bn
\ub{Stage B}:   We assume $(*)$. \newline
Let  $C = \{ \delta < \chi:\delta  = \text{ sup}(B \cap \delta)
\text{ and }  (s \in [B \cap \delta]^n$ for some \newline

$\qquad \qquad \qquad n < 1 + n^* \Rightarrow M_s \cap \chi \subseteq 
\delta) \}$. \newline
Let  $\gamma(*) = \text{ Min}(B)$.  Now for $p \in  P^*_\chi \cap 
M_{\{\gamma(\ast)\}}$ and $\bar c = \langle c_1,c_2 \rangle \in \sigma \times
\sigma$ let us define the statement 
\medskip
\roster
\item "{$(*)^{\bar c}_p$}"  if  
$p \le p^0 \in P^\ast \cap M_{\{ \gamma(*) \}}$
then we can find  $p^1,p^2 \in P^\ast_\chi \cap
M_{\{ \gamma(*)\}},p^0 \le p^1$,  $p^0 \le p^2$ such that for $\ell =1,2$:    
\newline
for  $\gamma_1 < \gamma_2$,  $\gamma_1 \in B$, $\gamma_2 \in B$,
we can find  $r_1$, $r_2 \in  P^\ast \cap  M_{\{\gamma_1,\gamma_2\}}$
(so \newline
$\text{Dom}(r_\ell) \subseteq  M_{\{\gamma_1,\gamma_2\}} \cap \chi)$
such that:

$$
r_\ell \Vdash ``{\underset\tilde {}\to \tau}(\{\gamma_1,\gamma_2\}) = c_\ell"
$$
 
$$
r_\ell \restriction (\chi \cap M_{\{\gamma_\ell \}}) \le 
f_{\{\gamma(\ast)\},\{\gamma_\ell \}}(p^1) \text{ (for strong system:
equality) }
$$
 
$$
r_\ell \restriction (\chi \cap M_{\{\gamma_{3-\ell }\}}) \le
f_{\{\gamma(\ast)\},\{\gamma_{3-\ell }\}}(p^2) \text{ (for strong system:
equality) }.
$$
\endroster
\medskip

\noindent
As $|\sigma \times \sigma| < \mu$ and the relevant forcing notions are
$(< \mu)$-strategically complete, easily
${\Cal I} = \{ p \in P^* \cap M_{\{\gamma(\ast)\}}:\text{ for some } \bar c, 
(*)^{\bar c}_p$ hold$\}$ is a dense
subset of  $P^*_\chi \cap  M_{\{\gamma(\ast)\}}$,  but this partial forcing
satisfies the $\mu^+$-c.c.  Hence we can find ${\Cal I}^*  = \{ p_\zeta:
\zeta < \mu \} \subseteq {\Cal I}$,  a maximal antichain of  
$P^\ast_\chi \cap M_{\{\gamma(\ast)\}}$ hence of $P^\ast_\chi$ (as
${}^{\mu \ge}(M_{\{ \gamma(*)\}})$ is a subset of $M_{\{ \gamma(*)\}}$).
For $p \in {\Cal I}^*$ we can choose $c_1(p),c_2(p) \in \sigma$ such that: 
$(*)^{(c_1(p),c_2(p))}_p$ hold.
\bn
\ub{Stage C}:   As  $G_R$ was any subset of  $R$  generic over $V$ to which
$r^*$ belongs, there are $R$-names
${\underset\tilde {}\to \gamma}(*),\langle ({\underset\tilde {}\to p_\xi},
{\underset\tilde {}\to c_1}({\underset\tilde {}\to p_\xi}),
{\underset\tilde {}\to c_2}({\underset\tilde {}\to p_\xi})):
\xi < \mu \rangle$,  $\langle {\underset\tilde {}\to M_s}:s \in
[{\underset\tilde {}\to B}]^{< \aleph_0} \rangle$, \newline
$\langle {\underset\tilde {}\to f_{s,t}}:(s,t) \in 
\dsize \bigcup_{n < 1 + n^*}
\left( [{\underset\tilde {}\to B}]^n \times [{\underset\tilde {}\to B}]^n 
\right) \rangle$
forced by $r^*$ to be as above.  As $R$ is $\chi$-complete,
$\chi > 2^{\kappa + \lambda + 2^\mu}$, without loss of generality 
$r^\ast$ forces values \newline
$\gamma(\ast),M_\emptyset,M_{\{\gamma(\ast)\}}$,  
$\langle (p^\ast_\zeta,c^\ast_1(p^*_\zeta),c^\ast_2(p^*_\zeta)):
\zeta < \mu \rangle$. 

We now try to choose by induction on $\zeta \le \theta + 1$,  
$\bar Q^\zeta,\alpha^\zeta,\gamma^\zeta$ such that:
\mr
\widestnumber\item{$(A)(a)$}
\item "{$(A)(a)$}"  $\bar Q^\zeta  \in  R$
\sn
$\qquad$ \item "{$(b)$}"  $\bar Q^0 = \{r^*\}$
\sn
$\qquad$ \item "{$(c)$}"  $\ell g(\bar Q^\zeta)  = \alpha^\zeta$ 
\sn
$\qquad$ \item "{(d)}"  $\xi  < \zeta  \Rightarrow  \bar Q^\xi  = 
\bar Q^\zeta \restriction \alpha^\xi$ 
\sn
$\qquad$ \item "{(e)}"  $\langle \alpha^\zeta:\zeta \le \theta + 1 \rangle$
is (strictly) increasing continuous
\sn
$\qquad$ \item "{(f)}"  $\alpha^\zeta < \gamma_\zeta < \alpha^{\zeta +1}$ 
\sn
$\qquad$ \item "{(g)}"  $\bar Q^{\zeta +1} \Vdash_R ``\gamma^\zeta  \in  
{\underset\tilde {}\to B}$"
\sn
$\qquad$ \item "{(h)}"  $\bar Q^{\zeta +1}$ 
forces $(\Vdash_R)$ a value to \newline
$\langle 
{\underset\tilde {}\to M_s} \cap V:s \in [{\underset\tilde {}\to B} \cap
(\gamma_\zeta + 1)]^{<1+n^*} \rangle$ which we call $\langle M_s:s \in 
[B_\zeta]^{< 1+n^*} \rangle$.
\ermn
$(B)$  if $\zeta \le \theta + 1$, cf$(\zeta) > \mu$ then:
\roster
\item "{$(a)$}"  $a^{\bar Q^{\zeta +1}}_{\alpha_\zeta}  = 
\bigcup \{ \chi \cap  
M_{\{\xi_1,\xi_2\}}:\{\xi_1,\xi_2\} \in  [\{\gamma_\epsilon:\epsilon  < 
\zeta \}]^{< 1+n^*} \}$
\sn
\item "{$(b)$}"  $I^{\bar Q^{\zeta + 1}}_{\alpha_\zeta} = \{b:b
\text{ an initial segment of } a^{\bar Q^{\zeta +1}}_{\alpha_\zeta} 
\text{ and } \text{cf(otp}(b)) \ne \mu^+\}$
{\roster
\itemitem{ {{}} }  [explanation: this satisfies the simplicity demands]
\endroster}
\item "{$(c)$}"  ${\underset\tilde {}\to Q^{\bar Q^{\zeta +1}}_{\alpha_\zeta}}
 = \{ h:h \text{ a function},
\text{Dom}(h) \subseteq \mu,|\text{Dom}(h)| < \mu,h(i) \in
Q^{\bar Q^\zeta}_{\alpha_\zeta,*}$ \newline

$\qquad \qquad \qquad \quad$ when defined $\}$ (see $(d)$ below) \newline
order $h_1 \le h_2$ if 
$i \in \text{ Dom}(h_1) \Rightarrow h_1(i) \subseteq h_2(i)$ where \nl
$Q^{{\bar Q}^\zeta}_{\alpha_\zeta,*}$ is defined in clause (d) below 
{\roster
\itemitem{ {{}} }  [explanation: the forcing notion in clause 
(d) adds a subset $\underset\tilde {}\to u$ of $\zeta$ such that on
$\{\gamma_\zeta:\zeta \in \underset\tilde {}\to u\}$ the 
colouring $\underset\tilde {}\to \tau$ get only two values; the 
forcing notion from clause (c) makes $\zeta$ the union of $\le \mu$ 
such sets and this induces a representation of $B_\zeta$ as a union of 
$\mu$ sets on each $\underset\tilde {}\to \tau$ get at most two colours]
\endroster}
\sn
\item "{$(d)$}"  $Q^{\bar Q^\zeta}_{\alpha_\zeta,*} = \biggl\{ 
u:u \in [\zeta]^{< \mu}$, and for some
$\xi < \mu$ we have: \newline

$\qquad \qquad \qquad \quad$ for every $j_1 < j_2$ from $u$,
we can find  $p^1$, $p^2$, $r_1$, $r_2$ \newline

$\qquad \qquad \qquad \quad$ such that for $\ell = 1,2$ we have:

$$
\qquad p^*_\xi \le  p^\ell \in M_{\{\gamma(*)\}} \cap P^*_\chi,        
$$

$$
\qquad r_\ell  \in  P^\ast _\chi  \cap  
M_{\{\gamma_{j_1},\gamma_{j_2}\}},
$$        

$$
r_\ell \Vdash ``{\underset\sim {}\to \tau} (\{\gamma_{j_1},\gamma_{j_2}\}) =
c^*_\ell(p^*_\xi),
$$

$$
\qquad \qquad \qquad r_\ell 
\restriction (\chi \cap  M_{\{\gamma_{j_\ell }\}}) \le 
f_{\{\gamma(*)\},\{\gamma_{j_1}\}}(p^1),
$$        

$$
\qquad \qquad \qquad
r_\ell \restriction (\chi \cap  M_{\{\gamma_{j_{3-\ell }}\}}) \le 
f_{\{\gamma(\ast)\},\{\gamma_{j_{3-\ell }}\}}(p^2),
$$        

$$
\qquad \qquad \text{and }
r_1 \in {\underset\tilde {}\to G_{P_{\alpha_\zeta}}} \text{ or }  r_2 \in  
{\underset\tilde {}\to G_{P_{\alpha_\zeta}}} \biggr\}.
$$
\endroster
\bn
\ub{Stage D}:  Again we shall use less than obtained for later use.

The point is to verify that we can carry the induction.  Now there is no
problem to do this for $\zeta=0$ and for $\zeta$ limit.  So we deal with
$\zeta +1,\zeta \le \theta$ and we are assuming that $\bar Q^\zeta$ is already
defined.  If cf$(\zeta) \le \mu$ clause (B) is empty and it is easy to
satisfy clause (A) is easy.  So assume cf$(\zeta) \ge \mu^+$.  Now as before
clause (A) is easy.  The point is to choose $\bar Q^{\zeta +1}$ or just
$\bar Q^{\zeta +1} \restriction (\alpha_\zeta +1)$ to satisfy clause (B).
Now ${\underset\tilde {}\to Q_{\alpha_\zeta}}$ is chosen by clause (B) so
$\bar Q^{\zeta +1} \restriction (\alpha_\zeta +1)$ is now fixed. \nl
The point is to prove that the condition concerning 
$\ast^\epsilon_\mu$ from Definition \scite{1.2} holds as required in
Definition \scite{1.7}(1)(d).  From now on we may omit the superscript 
$\bar Q^{\zeta + 1}$ or $\bar Q^{\zeta +1} \restriction (\alpha_\zeta +1)$ so
$P^*_{\alpha_\zeta} = P^{\bar Q^{\zeta + 1}}_{\alpha_\zeta} \restriction
(\alpha_\zeta +1)$, etc.

That is, we assume $b \in I_{\alpha_\zeta}$ and we will prove that
$(P^*_b,P^*_{a_{\alpha_\zeta} \cup \{ \alpha_\zeta\}})$ satisfies
$*^\varepsilon_\mu$. \newline
Note
\medskip
\roster
\item "{$(*)_1$}"  if $\bar Q^{\xi +1}$ is well defined (or just
$\bar Q^{\xi + 1} \restriction (\alpha_\xi +1) \in R)$ and
cf$(\xi) > \mu$ then ($P_{\alpha_\xi}$ is well defined and) in 
$V^{P_{\alpha_{\xi + 1}}},\{ \gamma_\Upsilon:\Upsilon < \xi \}$ is well
defined and it can be represented as $\dsize \bigcup_{i < \mu}
{\Cal U}_i$, such that each $u \in [{\Cal U}_i]^{< \mu}$ belongs to
$Q^{\bar Q^\xi}_{\alpha_\xi,*}$
\sn
\item "{$(*)_2$}"  if $\zeta(1) < \zeta(2) \le \zeta$ and cf$(\zeta(1)),
cf(\zeta(2)) > \mu$ then $Q^{\bar Q^{\zeta(1)}}_{\alpha_{\zeta(1)},\ast} 
\subseteq Q^{\bar Q^{\zeta(2)}}_{\alpha_{\zeta(2)},\ast}$, also for the 
compatibility relation
\sn
\item "{$(*)_3$}"  the elements of $Q^{\bar Q^\zeta}_{\alpha_\zeta,*}$
are from $V$, in fact are sets of ordinals of cardinality $< \mu$ ordered by
$\subseteq$ and the lub of set of cardinality $< \mu$ members is the union 
(if there is an upper bound), so 
$Q^{\bar Q^{\zeta}}_{\alpha_\zeta,\ast}$ is $\mu$-complete
\sn
\item "{$(*)_4$}"  $\bar Q^\zeta$ is well defined \underbar{and} \newline
$\Vdash_{P_{\alpha_\zeta}}$ ``for $\xi < \zeta$, if cf$(\xi) > \mu$ then,
\newline

$\qquad Q^{\bar Q^\xi}_{\alpha_\xi,*}$ is the union of $\mu$ sets,
each set $(< \mu)$-directed \newline

$\qquad$ and with any two elements having a lub".
\endroster
\medskip

\noindent
Hence
\medskip
\roster
\item "{$(*)_5$}"  if cf$(\zeta) > \mu^+$, then in $V^{P_{\alpha_\zeta}}$,
each subset of $Q^{\bar Q^{\zeta + 1}}_{\alpha_\zeta,*}$ of cardinality 
$\le \mu^+$ is included in the union of $\mu$
sets, each directed and with any two elements having a lub.
\endroster
\medskip
\noindent
Note that by the definition of $Q^{\bar Q^\zeta}_{\alpha_\zeta,*}$ we have
\mr
\item "{$(*)_6$}"  a family of $< \mu$ members of
$Q^{\bar Q^\zeta}_{\alpha_\zeta,*}$ has a common upper bound iff any two of
them are compatible, and then the union is a lub of the family.
\ermn
So if cf$(\zeta) > \mu^+$, we are done as by $(*)_5 + (*)_6$ we have
$\Vdash_{P_\zeta} ``Q_\zeta$ satisfies $*^\varepsilon_\mu"$ and can use
\scite{1.3}(4).
\medskip

So we can assume $\zeta = \Upsilon(*) \le \theta + 1$ and cf$(\zeta) =
\text{ cf}(\alpha_\zeta) = \mu^+$, and let \newline
$\langle \Upsilon(i):i < \mu^+ \rangle$ be increasing continuous with limit 
$\zeta$ and cf$(\Upsilon(i)) \le \mu$ for $i < \mu^+$. 
Let $b \in I_{\alpha_\zeta}$, hence $b$ is a bounded subset of $a_\zeta$. 
So by the induction hypothesis and \scite{1.3}(4) without loss of generality
$b = \bigcup\{ M_{\{\gamma_{\Upsilon_0},\gamma_{\Upsilon_1}\}} \cap 
\alpha_\zeta:\Upsilon_0 < \Upsilon_1 < \Upsilon(0)\}$.  \newline
\smallskip
\noindent
Define $c_0 = b_0 = b$ and for $\Upsilon \in [\Upsilon(0),\Upsilon(*))$ let
\newline
$b_{1,\Upsilon} = b_0 \cup \left( M_{\{\gamma_\Upsilon\}} \cap \alpha_\zeta 
\right) \cup \dsize \bigcup_{\Upsilon_1 < \Upsilon(0)} 
\left( M_{\{ \gamma_{\Upsilon_1},\gamma_\Upsilon\}} \cap \alpha_\zeta \right)$
\newline
(the third term could
be waived with minor changes), \newline
\smallskip
\noindent
$b_1 = b_{1,\Upsilon(0)}, b_2 = b_1 \cup b_{1,\Upsilon(0)+1},
c_2 = \bigcup\{b_{1,\Upsilon}:\Upsilon \in [\Upsilon(0),\zeta)\}$
\newline
\smallskip
\noindent
$c_3 = a_{\alpha_{\Upsilon(*)}} = \bigcup \{
M_{\{ \gamma_{\Upsilon_1},\gamma_{\Upsilon_2}\}} \cap \alpha_{\Upsilon(*)}:
\Upsilon_1 < \zeta,\Upsilon_2 < \zeta\}$ \newline
and $c_4 = a_{\alpha_{\Upsilon(*)}} \cup \{ \alpha_\zeta\}$.
\medskip
\noindent
\underbar{Note}:  There is no $c_1$. \newline 
All these sets are $\bar Q^{\alpha_\zeta + 1}$-closed.  We now choose
several winning strategies which exist by the induction hypothesis on
$\zeta$.
\medskip

Let $St_0$ be a winning 
strategy of the first player in a game above $*^\varepsilon_\mu
[P^*_{b_0}]$.  Let $St_1$ be a winning strategy of the first player in
$*^\varepsilon_\mu[P^*_{b_0},P^*_{b_1}]$ which projects to $St_0$.   
For every $\Upsilon \in
[\Upsilon(0),\Upsilon(*))$ let $St_{1,\Upsilon}$ be a winning strategy of 
the first player in $*^\varepsilon_\mu[P^*_{b_0},P^*_{b_{1,\Upsilon}}]$ 
conjugate to $St_1$ (by OP$_{b_{1,\Upsilon},b_1})$.
\medskip

For $\bar \Upsilon = \langle \Upsilon_1,\Upsilon_2 \rangle,\Upsilon_1
< \Upsilon_2,\{ \Upsilon_1,\Upsilon_2\} \subseteq [\Upsilon(0),\Upsilon(*))$ 
let \newline
$b_{2,\bar \Upsilon} = b_{1,\Upsilon_1} \cup b_{1,\Upsilon_2} \cup
\left( M_{\{ \Upsilon_1,\Upsilon_2\}} \cap \alpha_\zeta \right)$ and
let $St_{2,\bar \Upsilon}$ be a winning strategy in $*^\varepsilon_\mu
[P^*_{b_{1,\Upsilon_1} \cup b_{1,\Upsilon_2}},P^*_{b_{2,\bar \Upsilon}}]$ 
which is above $St_{1,\Upsilon_1} \times St_{1,\Upsilon_2}$ (remember that 
both project to $St_0$); also note
as long as the second player uses conditions in 
$P^*_{b_{\ell,\Upsilon_\ell}}$ then so does the first player (for each
$i < \mu^+$ separately).
\newline
Also, the first player has a winning strategy in $*^\varepsilon_\mu[P^*_{c_0},
P^*_{c_2}]$ but we want a very special winning strategy $St_2$:
(letting $g_2$ be a fixed pairing function on $\mu^+$) in a play
$\bigl< \langle p^\xi_i:i < \mu^+ \rangle,\langle q^\xi_i:i < \mu^+ \rangle,
f^\xi:\xi < \varepsilon \bigr>$ where the first player uses the strategy
$St_2$ we demand that clauses $(a)-(d)$ below holds:
\medskip
\roster
\item "{$(a)$}"  $\bigl< \langle p^\xi_i \restriction b_0:i < \mu^+ \rangle,
\langle q^\xi_i \restriction b_0:i < \mu^+ \rangle,f^{1,\xi}:\xi <
\varepsilon \bigr>$ \nl
is a play of $*^\varepsilon_\mu[P^*_{b_0}]$ in which the
first player uses the strategy $St_0$)
\smallskip
\noindent
\item "{$(b)$}"  for each $\Upsilon \in [\Upsilon(0),\Upsilon(*))$ defining
$$
p^{2,\Upsilon,\xi}_i = \cases p^\xi_i \restriction b_{1,\Upsilon}
&\text{ \underbar{if} } \quad \Upsilon(i) > \Upsilon \\
  p^\xi_i \restriction b_0 &\text{ \underbar{if} } \quad 
\Upsilon(i) \le \Upsilon
\endcases
$$
\medskip
$$ 
q^{2,\Upsilon,\xi}_i = \cases q^\xi_i \restriction b_{1,\Upsilon}
&\text{ \underbar{if} } \quad \Upsilon(i) > \Upsilon \\
  q^\xi_i \restriction b_0 &\text{ \underbar{if} } \quad 
\Upsilon(i) \le \Upsilon
\endcases
$$
\endroster
\medskip

\noindent
we have: $\bigl< \langle p^{2,\Upsilon,\xi}_i:i < \mu^+ \rangle,
\langle q^{2,\Upsilon,\xi}_i:i < \mu^+ \rangle,f^{2,\Upsilon,\xi}:
\xi < \varepsilon \bigr>$ is a play of $*^\varepsilon_\mu
[P^*_{b_0},P^*_{b_1,\Upsilon}]$ in which the first player uses the strategy 
$St_{1,\Upsilon}$.
\medskip
\roster
\item "{$(c)$}"  For any pair $\bar \zeta = (\zeta_1,\zeta_2)$ of ordinals
in $\mu \times \varepsilon$, let

$$
\Upsilon(i,\bar \zeta) = \Upsilon_{\bar \zeta}(i) \text{ is the } 
\zeta_1 \text{-th member of Dom}(q^{\zeta_2}_i) \backslash \Upsilon(i)
$$

$$
p^{3,\bar \zeta,\xi}_i = \text{ OP}_{b_{1,\Upsilon(0)},b_{1,\Upsilon_{\bar
\zeta}(i)}}(p^\xi_i \restriction b_{1,\Upsilon_{\bar \zeta}(i)})
$$

$$
q^{3,\bar \zeta,\xi}_i = \text{ OP}_{b_{1,\Upsilon(0)},b_{1,\Upsilon_{\bar
\zeta}(i)}}(q^\xi_i \restriction b_{1,\Upsilon_{\bar \zeta}(i)}),
$$
\endroster
\medskip

\noindent
we demand that
$\bigl< \langle p^{3,{\bar \zeta},\xi}_i:i < \mu^+ \rangle,
\langle q^{3,{\bar \zeta},\xi}_i:i < \mu^+ \rangle,f^{3,{\bar \zeta},\xi}:
\xi < \varepsilon \bigr>$ is a play of \newline
$*^\varepsilon_\mu[P^*_{b_0},P^*_{b_{1,\Upsilon(0)}}]$ in which
the first player uses the strategy $St_{1,\Upsilon(0)}$.

So for each $i < \mu$, for $\zeta_1 < \mu$ too large 
$\Upsilon(i,\bar \zeta)$ is not well defined
and we stipulate the forcing conditions are $\emptyset$.
\medskip
\roster
\item "{$(d)$}"  $f^\xi(i)$ codes $f^{1,\xi}(i),
\langle f^{2,\Upsilon,\xi}(i):\Upsilon \in [\Upsilon(0),\Upsilon(*))$ and
$(\exists \beta \in b_{1,\Upsilon}
\backslash b_0)[p^\xi_i(\beta) \ne 0_{\underset\tilde {}\to Q_\beta}]\rangle$
and $\langle f^{3,\bar \zeta,\xi}(i):\bar \zeta \in \mu \times \varepsilon,
\text{ and } \Upsilon_{\bar \zeta}(i) \text{ is well defined}\rangle$ and
the information on $p^\xi_i(\alpha_{\Upsilon(*)})$ and it codes
$$
\align
\bigl\{\langle j_1,\zeta_1,\zeta_2 \rangle:&\beta,\text{ the } 
\zeta_2 \text{-th member of Dom}(p^\xi_i) \\
  &\text{satisfies}: j_1 = \text{ Min}\{j:\beta \in \text{ Dom}(p^\xi_j)\}, \\
  &\beta \text{ is the } \zeta_1 \text{-th member of Dom}(p^\xi_{j_1}) 
\bigr\}
\endalign
$$
\medskip

\noindent
and

$$
\align
\bigl\{ \langle j,\zeta_1,\zeta_2 \rangle:&\text{ for some } \Upsilon,
\beta, \text{ the } \zeta_1 \text{-th member of Dom}(p^\xi_i), \\
  &\text{ belongs to } b_{1,\Upsilon} \backslash b_0 \text{ and satisfies}: \\
  &\,j = \text{ Min}\{j':(\text{Dom}(p^\xi_{j'}) \cap 
(b_{1,\Upsilon} \backslash b_0) \ne \emptyset\} \\
  &\text{ and the } \zeta_2 \text{-th member of Dom}(p^\xi_j)
\text{ belongs to } b_{1,\Upsilon} \backslash b_0 \bigr\}
\endalign
$$
\endroster
\medskip

\noindent
(note: for each $\zeta_2 < \varepsilon,i < \mu^+$ we have: \newline
$\{ \zeta_1 < \mu:\Upsilon_{(\zeta_1,\zeta_2)}(i) \text{ is well defined}\}$ 
is a bounded subset of $\mu$).

Check that such $St_2$ exists, (note that the number of times we have 
to increase $p_i \restriction b_0$ is $< \mu$).
\medskip

Clearly $c_2 \subseteq c_3$ are $\bar Q$-closed, hence there is a winning
strategy $St_3$ of the first player in 
$*^\varepsilon_\mu[P^*_{c_2},P^*_{c_3}]$ above $St_2$.
\medskip
\roster
\item "{$(e)$}"  For any 
$\bar \Upsilon = (\Upsilon_1,\Upsilon_2)$ such that $\Upsilon(0) \le
\Upsilon_1 < \Upsilon_2 < \Upsilon(*)$, and defining \newline
$p^{4,\bar \Upsilon,\xi} = p^\xi_i \restriction b_{2,\bar \Upsilon}$, \newline
$q^{4,\bar \Upsilon,\xi}_i = q^\xi_i \restriction b_{2,\bar \Upsilon}$
\newline
(can behave similarly in clause $(b)$), \nl
\ub{we have}: $\bigl< \langle p^{4,\bar \Upsilon,\xi}_i:i < \mu^+ \rangle,
\langle q^{4,\bar \Upsilon,\xi}_i:i < \mu^+ \rangle,f^{4,{\bar \Upsilon},\xi}:
\xi < \varepsilon \bigr>$ \nl
is a play of $*^\varepsilon_\mu[P^*_{b_{2,\bar \Upsilon}}]$ 
in which the first player uses the strategy $St_{2,\bar \Upsilon}$.
\ermn
Lastly, let $St_4$ be a strategy of the first player in $*^\varepsilon_\mu
[P^*_{c_3},P^*_{c_4}]$ which is above $St_3$ and it guarantees:
\medskip
\roster
\item "{$(*)$}"  if $\bigl< \langle p^\xi_i:i < \mu^+ \rangle,
\langle q^\xi_i:i < \mu^+ \rangle,f^4_\xi:\xi < \varepsilon \bigr>$ is a
play of the game in which the first player uses his strategy $St_4$ 
\underbar{then}:
{\roster
\itemitem{ $(\alpha)$ } $q^\xi_i \restriction a_{\alpha_\Upsilon}$ forces
a value to $p^\xi_i(\alpha_{\Upsilon(*)})$
\itemitem{ $(\beta)$ }  if $\Upsilon_1 \ne \Upsilon_2$ are 
from (the value forced on) $q^\xi_i(\alpha_{\Upsilon(*)})$ \underbar{then} 
$q^\xi_i \restriction a_\Upsilon$ is above the relevant parts of witnesses to
this.
\endroster}
\endroster
\medskip

Clearly $St_4$ is (essentially) a strategy of the first player in
$*^\varepsilon_\mu[P^*_{b_0},P^*_{c_4}]$ (for the almost
$*^\varepsilon_\mu$ case above $St_0$).  All we have to prove
is that $St_4$ is a winning strategy.  So let
$\bigl< \langle p^\xi_i:i < \mu^+ \rangle$, 
$\langle q^\xi_i:i < \mu^+ \rangle,f^4_\xi:\xi < \varepsilon \bigr>$ be a
play of $*^\varepsilon_\mu[P^*_{b_0},P^*_{c_4}]$ in which the first player 
uses the strategy $St_4$. \newline
By the definition of the game $*^\varepsilon_\mu
[P^*_{b_0},P^*_{c_4}]$ \wilog \, for some club $E_1$ of $\mu^+$ 
(see clause $(a)$):
\medskip
\roster
\item "{$(**)_1$}"  if $\{i,j\} \subseteq S^{\mu^+}_\mu \cap E_1$ 
(see \sciteu{0.I}) and \newline
$\dsize \bigwedge_{\xi < \varepsilon} f^4_\xi(i) = f^4_\xi(j)$ \underbar{then}
\newline
$\{ p^\xi_i \restriction b_0,p^\xi_j \restriction b_0:\xi < \varepsilon\}$
has an upper bound in $P^*_{b_0}$.
\endroster
\medskip

\noindent
By clause $(b)$ in the demands on $St_{1,\Upsilon}$ for some club $E_2$ of 
$\mu^+$ we have:
\medskip
\roster
\item "{$(**)_2$}"  if $\{i,j\} \subseteq S^{\mu^+}_\mu \cap E_2$ and \newline
$\Upsilon \in [\Upsilon(0),\Upsilon(*))$ and \newline
$\dsize \bigwedge_{\xi < \varepsilon}[(b_{1,\Upsilon} \backslash b_0) \cap
\text{ Dom}(p^\xi_i) \ne \emptyset \and (b_{1,\Upsilon} \backslash b_0) \cap
\text{ Dom}(p^\xi_j) \ne \emptyset \Rightarrow 
f^{2,\Upsilon,\xi}(i) = f^{2,\Upsilon,\xi}(j)]$ \newline
(which holds if
$\dsize \bigwedge_{\xi < \varepsilon}
f^4_\xi(i) = f^4_\xi(j)$), and $r$ is an upper bound of \newline
$\{ p^\xi_i \restriction b_0,p^\xi_j \restriction b_0:\xi < \varepsilon\}$
\underbar{then}
\newline
$\{ p^\xi_i \restriction b_{1,\Upsilon},p^\xi_j \restriction b_{1,\Upsilon}:
\xi < \varepsilon\} \cup \{r\}$ has an upper bound in $P^*_{b_{1,\Upsilon}}$.
\endroster
\medskip

\noindent
By clause $(c)$ in the choice of $St_2$ we know that there is a club $E_3$
of $\mu^+$ such that:
\medskip
\roster
\item "{$(**)_3$}"  if $\bar \zeta \in \mu \times \varepsilon,
\{i,j\} \subseteq S^{\mu^+}_\mu \cap E_3$ and \newline
$\dsize \bigwedge_{\xi < \varepsilon} f^{3,\bar \zeta,\xi}(i) = 
f^{3,\bar \zeta,\xi}(j)$ (which holds if $\dsize \bigwedge_{\xi < \varepsilon}
f^4_\xi(i) = f^4_\xi(j)$) and \newline
$r \in P^*_{b_0}$ is an upper bound of 
$\{ p^\xi_i \restriction b_0,q^\xi_i \restriction b_0:\xi < \zeta\}$
\underbar{then}
\newline
$\{ p^{3,\bar \zeta,\xi}_i,p^{3,\bar \zeta,\xi}_j:
\xi < \varepsilon\} \cup \{r\}$ has an upper bound.
\endroster
\medskip

\noindent
By clause $(e)$ in the demand on $St_3$, for some club $E_4$ of $\mu^+$
\medskip
\roster
\item "{$(**)_4$}"  if $\{i,j\} \subseteq S^{\mu^+}_\mu \cap E_4$ and
$\dsize \bigwedge_{\xi < \varepsilon} f^4_\xi(i) = f^4_\xi(j)$ and $r$ is an
upper bound of $\{p^\xi_i \restriction b_0,p^\xi_j \restriction b_0:\xi <
\varepsilon\}$ then $\{p^\xi_i \restriction \Upsilon(i),p^\xi_j \restriction
\Upsilon(j):\xi < \varepsilon\} \cup \{r\}$ has an upper bound in
$p^*_{\alpha_0}$ (even $p^*_{\alpha_{\text{Max}\{\Upsilon(i),
\Upsilon(j)\}}}$).
\endroster
\medskip

\noindent
Last
\medskip
\roster
\item "{$(**)_5$}"  $E$ is a club of $\mu^+$ included in $E_1 \cap E_2 \cap
E_3 \cap E_4$ such that: \newline
$i < j \in E \Rightarrow \text{ Dom}(p^\xi_i \restriction
c_3) \cup \text{ Dom}(q^\xi_i \restriction c_3) \subseteq
\alpha_{\Upsilon(j)}$.
\endroster
\bigskip

\noindent
The rest is as in \cite[\S2]{Sh:276}. \hfill$\square_{1.12}$
\enddemo
\bigskip

\proclaim{\stag{1.13} Theorem}  We can in \scite{1.12} 
replace ``measurable", by (strongly) Mahlo.
\endproclaim
\bigskip

\remark{\stag{1.13A} Remark}  It is not straightforward; e.g. we may use the
version of squared diamond given in Fact \scite{1.15} below.

We first prove two claims. 
\endremark
\bigskip

\proclaim{\stag{1.14} Claim}  Suppose $\lambda$ is a strongly 
inaccessible Mahlo cardinal, $\chi > \lambda > \theta = \theta^{< \sigma}$,
${\frak C}$ an expansion of $({\Cal H}(\chi),\in,<^*_\chi)$ by $\le \theta$
relations.  \underbar{Then} for some club $E$ of $\lambda$ for every 
inaccessible $\kappa \in E$ we have:
\medskip
\roster
\item "{$(*)_\kappa$}"  for every $x \in {\Cal H}(\chi)$ there are 
$B \in [\kappa]^\kappa$ and $N_s$
(for $s \in [B \cup \{ \kappa \}]^2),N'_{\{ i \}}$ \newline
(for $i \in B \cup \{ \kappa \}),N_{\{ i \}}$ 
(for $i \in B)$ and $N_\emptyset$ (so 
$N_{\{ \kappa \}}$ is meaningless) such that ($L_{\sigma,\sigma}$ is like
the first order logic but with conjunctions and a string of existential
quantifiers of any length $< \sigma$):
{\roster
\itemitem{ (a) }  $x \in N_s \prec_{L_{\sigma,\sigma}} {\frak C}$ and
$\theta \subseteq N_s$
sn
\itemitem{ (b) }  $x \in N'_{\{ i \}} \prec_{L_{\sigma,\sigma}}
{\frak C}$ and $\theta \subseteq N'_{\{ i \}} \subseteq N_{\{ i \}}$
\sn
\itemitem{ (c) }  $s \subseteq B \Rightarrow N_s \cap \lambda \subseteq
\kappa \and N'_s \cap \lambda \subseteq \kappa$ (when defined)
\sn
\itemitem{ (d) }  $N_\emptyset \prec_{L_{\sigma,\sigma}} N_{\{ i \}}$ and  
$\text{min}(N_{\{ i \}} \cap \lambda \backslash N_\emptyset) > 
\text{ sup}\left[ \bigcup \{ N_s \cap \lambda:s 
\subseteq [B \cap i]^{\le 2} \} \right]$
\sn
\itemitem{ (e) }  for $j < i,N_{\{ j,i \}}$ is the $L_{\sigma,\sigma}$-
Skolem hull of $N_{\{ j \}} \cup N'_{\{ i \}}$ inside ${\frak C}$
\sn
\itemitem{ (f) }  for $j < i,N_{\{ j \}} \cap \lambda$ is an
initial segment of $N_{\{ j,i \}} \cap \lambda$
\sn
\itemitem{ (g) }  for $j < i$, Min$\left( N_{\{j,i\}} \cap \lambda \backslash
N_{\{j\}} \right) > \sup \{N_{\{j_1,i_1\}} \cap \lambda:j_1 < i_1 < i\}$
\sn
\itemitem{ (h) }  $N_s,N'_s$ have cardinality $\theta$ when defined.
\endroster}
\endroster
\endproclaim
\bigskip

\demo{Proof}  Let $\theta_1 = 2^\theta,\theta_2 = 2^{\theta_1}$.
Let ${\frak A}$ and $\kappa$ be such that:
$$
\kappa \text{ strongly inaccessible}
$$

$$
{\frak A} \prec_{L_{\theta^+_2,\theta_2}} {\frak C}
$$

$$
{\frak A}^{< \kappa} \subseteq {\frak A}
$$

$$
{\frak A} \cap \lambda = \kappa.
$$

\noindent
(Clearly for some club $E$ of $\lambda$, for every strongly inaccessible
$\kappa \in E$ there is ${\frak A}$ as above; so it is enough to prove
$(*)_\kappa$).  Without loss of generality, $\kappa > \theta$.
Next choose ${\frak B}_i \prec_{L_{\theta^+_2,\theta^+_2}} {\frak C}$, 
increasing continuous in $i$ for $i < \kappa,\langle {\frak B}_i:i \le 
j \rangle \in {\frak B}_{j + 1},\| {\frak B}_j \| < \kappa,{\frak B}_i 
\cap \kappa$ an ordinal and $\{x,\lambda,\theta,\sigma,\kappa,\lambda,
{\frak A}\} \in {\frak B}_0$. \newline
Let ${\frak B} = {\frak B}_{\theta^+}$, and let $f$ be a function from 
${\frak B}$ into ${\frak A}$, which is an $\prec_{L_{\theta^+_1,\theta^+_2}}$
elementary mapping (for the model 
${\frak C},\text{Dom}(f) = {\frak B},\text{Rang}(f) \subseteq {\frak A}$).
\medskip

Let $N \prec_{L_{\sigma,\sigma}} {\frak C}$ be such that \newline
$\{x,{\frak A},{\frak B},\langle {\frak B}_i:i \le \theta^+ \rangle,
f,\sigma,\theta,\lambda,\kappa \} \in N,\theta + 1 
\subseteq N,\| N \| = \theta,N^{< \sigma} \subseteq N$.

Let $N^+$ be the $L_{\sigma,\sigma}$-Skolem hull of $N \cup f(N)$ in
${\frak C}$. \newline
Let $N_\emptyset$ be $N^+ \cap {\frak A} \cap {\frak B}$, as
$\| N_\emptyset \| \le \theta$ we have $N_\emptyset \in {\frak A} \cap
{\frak B}$.  Let $N_{\{ 0 \}} = N^+ \cap {\frak A}$ (so
$N_\emptyset = N_{\{ 0 \}} \cap {\frak B}$, and $N_\emptyset \cap \lambda
(\subseteq \kappa)$ is an initial segment of $N_{\{ 0 \}} \cap \lambda
(\subseteq \kappa)$, let $N'_{\{ \kappa \}} = N^+ \cap {\frak B}$ and
$N'_{\{ 0 \}} = f(N'_{\{ \kappa \}})$, so $N'_{\{ 0 \}} \prec N_{\{ 0 \}}$.
Let $\alpha_0 = f(\kappa)$.  Now we choose by induction on $i < \kappa,
\alpha_i,N'_{\{ i \}},N_{\{ i \}},g_i$ and $N_{\{ i,j \}}$
for $j < i$ such that:
\medskip
\roster
\item  $g_i$ is an $\prec_{L_{\sigma,\sigma}}$ elementary mapping from 
$N_{\{ 0 \}}$ into ${\frak A},g_0 = \text{ id}_{N_{\{ 0 \}}}$
\sn
\item  $g_i(\alpha_0) = \alpha_i$
\sn
\item  for $j < i,N_{\{ j,i \}}$ is the $L_{\sigma,\sigma}$-Skolem hull of
$N_{\{ j \}} \cup N'_{\{ i \}}$ (in ${\frak C}$)
\sn
\item  $N_{\{ i,\kappa \}}$ is the $L_{\sigma,\sigma}$-Skolem hull of
$N_{\{ i \}} \cup N'_{\{ \kappa \}}$
\sn
\item  $N_{\{ i,\kappa \}},N_{\{ 0,\kappa \}}$ are isomorphic, in fact
there is an ismorphism from $N_{\{ 0,\kappa \}}$ onto $N_{\{ i,\kappa \}}$
extending $g_i \cup \text{ id}_{N'_{\{ \kappa \}}}$
\sn
\item  for $j < i$ there is an isomorphism from $N_{\{ j,i \}}$ onto
$N_{\{ j,\kappa \}}$ extending \newline
$\text{id}_{N_{\{ j \}}} \cup (f^{-1} \circ g^{-1}_i) \restriction 
N'_{\{ i \}}$
\sn
\item  $N_{\{ j \}} \cap \lambda$ is an initial segment of $N_{\{ j,i \}} \cap
\lambda$ for $j < i$.
\endroster
\medskip

\noindent
This is possible and gives the desired result.  \hfill$\square_{\scite{1.14}}$

\enddemo
\bigskip

\fakesubhead{\stag{1.15} Fact}  \endsubhead  Let 
$\chi$ be strongly inaccessible $(k+1)$-Mahlo, $\kappa < \chi$ are regular.  
By a forcing with a $P$ which is $\kappa^+$-complete of cardinality $\chi$, 
not collapsing cardinals nor cofinalities nor changing cardinal arithmetic 
we can get:
\medskip
\roster
\item "{$(*)^{\kappa,k}_\chi$}"  there is 
$\bar A = \langle A_\alpha:\alpha < \chi \rangle$ and $\bar C = \langle
C_\alpha:\alpha \in S \rangle$ such that:
{\roster
\itemitem{ (a) }  $S \subseteq \{ \delta < \chi:\delta > \kappa \text{ and }
\text{cf}(\delta) \le \kappa\}$ and 
$\{ \delta \in S:\text{otp}(C_\delta) = \kappa\}$ is a
stationary subset of $\chi$
\sn
\itemitem{ (b) }  $C_\alpha \subseteq \alpha \cap S,[\beta \in C_\alpha 
\Rightarrow C_\beta = C_\alpha \cap \beta],\text{ otp}(C_\alpha) \le \kappa,
C_\alpha$ a closed subset of $\alpha$ and $[\sup(C_\alpha) = \alpha
\Leftrightarrow C_\alpha$ has no last element]
\sn
\itemitem{ (c) }  $A_\alpha \subseteq \alpha$
\sn
\itemitem{ (d) }  $\beta \in C_\alpha \Rightarrow A_\beta = A_\alpha \cap
\beta$
\sn
\itemitem{ (e) }  $\bigl\{ \lambda < \chi: \lambda \text{ inaccessible}$, 
and for every $X \subseteq \lambda$ the set \newline

$\qquad \qquad \quad$ we have 
$\{ \alpha < \lambda:\text{otp}(C_\alpha) = \kappa,
X \cap \alpha = A_\alpha\}$ \newline

$\qquad \qquad \quad$ is a stationary subset of $\lambda \bigr\}$ \newline

$\,\,\,\,\,\,$is not only stationary but is a $k$-Mahlo subset, \nl

$\,\,\,\,\,\,$ moreover we actually get:
\sn
\itemitem{ (e)$^+$ }  for every strongly inaccessible $\lambda \in (\theta,
\chi),\langle (A_\alpha,C_\alpha):\alpha \in S \cap \lambda \rangle$ is
a club guessing squared diamond, that is clauses (a)-(d) hold with 
$\lambda,S \cap \lambda$ and: for every club $E$ of $\lambda$ and 
$X \subseteq \lambda$ for some $\delta \in S$ we have 
$C_\delta \cup \{ \delta\} \subseteq E$ and otp$(C_\delta) = \kappa$ and
$\alpha \in C_\delta \cup \{\delta\} \Rightarrow A_\alpha = X \cap \alpha$.
\endroster}
\endroster
\bigskip

\demo{Proof}  This can be obtained e.g. by iteration with Easton support,
in which for each strongly inaccessible $\lambda \in (\kappa,\chi]$ we 
add $\bar A,\bar C$ satisfying $(a)-(d)$ above, each condition being an 
initial segment.

More specifically, we define and prove by induction on $\alpha
\le \chi$

$$
\align
(1) \text{ [Definition}] \qquad P_\alpha = \biggl\{ (a,\bar C,\bar A):&(a) \,
a \subseteq \alpha \backslash \kappa^+, \\
  &(b) \, \text{ for every strongly inaccessible } \lambda \in 
          (\kappa,\chi] \\
  &\quad \text{ we have } \lambda > \sup(a \cap \lambda) \\
  &(c) \, \,\bar C = \langle C_\alpha:\alpha \in a \rangle \\
  &(d) \,C_\alpha \ne \emptyset \Rightarrow \text{ cf}(\alpha) \le \kappa
\and \text{ otp}(C_\alpha) \le \kappa \\
  &(e) \,\beta \in C_\alpha \Rightarrow \beta \in a \and 
C_\beta = C_\alpha \cap \beta \\
  &(f) \,C_\alpha \ne \emptyset \Rightarrow C_\alpha \text{ closed} \\
  &(g) \,\bar A = \langle {\underset\tilde {}\to A_\alpha}:\alpha \in a 
      \rangle \\
  &(h) \,{\underset\tilde {}\to A_\alpha} \text{ is a } 
P_\alpha \text{-name of a subset of } \alpha \\
  &(i) \,\beta \in C_\alpha \Rightarrow \Vdash_\alpha ``{\underset\tilde {}\to
A_\alpha} \cap \beta = {\underset\tilde {}\to A_\beta} \biggr\}
\endalign
$$
\smallskip

\noindent
\underbar{order} $p \le q$ iff $a^p \subseteq a^q,
\bar C^p = \bar C^q \restriction
a^p,\bar A^p = \bar A^q \restriction a^p$. \newline
\smallskip
\noindent
2)  [Claim]: $\beta < \alpha \Rightarrow P_\beta \lessdot P_\alpha$. \newline
3)  [Claim]: If $p \in P_\alpha,\beta < \alpha$, then $p \restriction \beta =
(a^p \cap \beta,\bar C \restriction (a \cap \beta),\bar A \restriction (a \cap
\beta))$ belongs to $P_\beta$ and: 
if $p \restriction \beta \le q \in P_\beta$ then $p,q$ are compatible 
in a simple way: $p \and q$ is a lub of $\{p,q\}$. \nl
4)  [Claim]: If $\lambda$ is strongly inaccessible $\le \chi$ and $> \kappa$
then $P_\lambda = \dbcu_{\alpha < \lambda} P_\alpha$.  If in addition 
$\lambda$ is Mahlo, then $P_\lambda$ satisfies the $\lambda$-c.c.

Let 
${\underset\tilde {}\to c_\alpha} = c^p_\alpha,
{\underset\tilde {}\to A_\alpha} = A^p_\alpha$ for every large enough 
$p \in {\underset\tilde {}\to G_{P_\chi}}$.
The point is that for every strongly inaccessible $\lambda \in (\theta,\chi],
P_\chi/P_\lambda$ does not add any subset of $\lambda$, and so $\langle
({\underset\tilde {}\to C_i},{\underset\tilde {}\to A_i}[G]):
i < \lambda \rangle$ is as required. \hfill$\square_{\scite{1.15}}$
\enddemo
\bigskip

\fakesubhead{\stag{1.16} Conclusion}  \endsubhead  Let 
$\theta = \theta^{< \sigma} < \lambda,\lambda$ a strongly inaccessible 
Mahlo cardinal, \ub{then} for some $\theta^+$-complete,
$\lambda$-c.c. forcing notion of cardinality $\lambda$ not collapsing
cardinals not changing cofinalities nor changing cardinal arithmetic, in
$V^P$ we get:
\medskip

\noindent
$(**)^{\theta,2}_\lambda$ $\,\,$ there are $\langle (B_\alpha,\bar M^\alpha,
C_\alpha):\alpha \in S \rangle$ such that:
\medskip
\roster
\item "{$(a)$}"  $S \subseteq \{ \delta < \chi:\text{cf}(\delta) \le \theta
\}$ and \newline
$\{ \delta \in S:\text{otp}(C_\delta) = \theta\}$ is a stationary subset of
$\chi$ and even of any strongly inaccessible $\lambda \in (\theta,\chi)$
\sn
\item "{$(b)$}"  $C_\alpha \subseteq \alpha \cap S,[\beta \in C_\alpha
\Rightarrow C_\beta = C_\alpha \cap \beta],\text{ otp}(C_\alpha) \le \theta,
C_\alpha$ a closed subset of $\alpha$ so $[\sup(C_\alpha) = \alpha
\Leftrightarrow C_\alpha$ has no last element)
\sn
\item "{$(c)$}"  $B_\alpha \subseteq \alpha,\text{ otp}(B_\alpha) = \omega
\times \text{ otp}(C_\alpha),\beta \in C_\alpha \Rightarrow B_\beta =
B_\alpha \cap \beta$
\sn
\item "{$(d)$}"  each $\langle M^\alpha_s:s \in [B_\alpha]^{\le 2} \rangle$ 
is as in \scite{1.14} (and $B_\alpha \subseteq B$) and \nl
$\beta \in C_\alpha \and s \in [B_\beta]^{\le 2} \Rightarrow M^\alpha_s =
M^\beta_s$
\sn
\item "{$(e)$}"  diamond property: if ${\frak B}$ is an expansion of
$({\Cal H}(\chi),\in,<^*_\chi)$ by $\le \theta$ relations, 
$B \in [\chi]^\chi$ then
for a club $E$ of $\chi$ for every strong inaccessible $\lambda \in 
\text{ acc}(E)$ for stationarily many $\delta \in S \cap \lambda$ we have
otp$(C_\delta) = \kappa,C_\delta \subseteq E$ and $B_\delta \subseteq B$ and
$s \in [B_\delta]^{\le 2} \Rightarrow M^\delta_s \prec {\frak B}$.
\endroster
\bigskip

\demo{Proof}  By \scite{1.15} + \scite{1.14} 
(alternatively, force this directly: simpler than in \scite{1.15}.
\enddemo
\bigskip

\remark{Remark}  In \scite{1.15} we could force a stronger version.
\endremark
\bigskip

\demo{Proof of \scite{1.13}}  We repeat the main proof the one of Theorem
\scite{1.12}, but using the diamond from \scite{1.15} for $k=0$.  In fact
the proof of \scite{1.12} was written such that it can be read as a proof of
\scite{1.13}, mainly in stage B we can get $(*)$ which is proved using
measurability, but use only $(*)'$.

\hfill$\square_{\scite{1.13}}$
\enddemo
\bigskip

\noindent
Combining the above proof and \cite{Sh:288} we get
\proclaim{\stag{1.17} Theorem}  Suppose
\medskip
\roster
\item "{$(a)$}"  $\mu = \aleph_0$ \underbar{or} $\mu$ is Laver indestructible
supercompact (see \cite{L}) or just $\mu$ as in \newline
\cite[\S4]{Sh:288}
\sn
\item "{$(b)$}"  $\lambda$ is $n^*$-Mahlo, $\lambda > \theta > \mu$
\sn
\item "{$(c)$}"  $k_{n^*}$ as in \cite{Sh:228} (see below).
\endroster
\medskip

\noindent
\underbar{Then} for some $\mu^+$-c.c. forcing notion $P$ of cardinality
$\lambda$ we have:

$$
\align
\Vdash_P ``&2^\mu = \lambda \rightarrow [\theta]^{n^*+1}_{k_{n^*}+1}",
\text{ moreover for } \sigma < \mu, \\
  &\lambda \rightarrow [\theta]^{n^*+1}_{\sigma,k_{n^*}}.
\endalign
$$

\endproclaim
\bigskip

\remark{\stag{1.17A} Remark}  1) What is $k_{n^*}$?
\mn
\ub{Case 1}:  $\mu = \aleph_0$; define on $[{}^\omega 2]^{n^*}$ an equivalence
relation $E$: if $w_1 = \{\eta_\ell:\ell < n^*\},w_2 = \{\nu_\ell:\ell <
n^*\}$ are members of $[{}^w 2]^{n^*}$ both listed in lexicographic 
increasing order, then $w_1 E w_2$ \ub{iff} for any $\ell_1 < \ell_2 < n^*$
and $\ell_3 < \ell_4 < n^*$ we have

$$
\ell g(\eta_{\ell_1} \cap \eta_{\ell_2}) < \ell g(\eta_{\ell_3} \cap
\eta_{\ell_4}) \Leftrightarrow \ell g(\nu_{\ell_1} \cap \nu_{\ell_2}) <
\ell g(\nu_{\ell_3} \cap \nu_{\ell_4}).
$$
\mn
Lastly, $k_{n^*}$ is the number of $E$-equivalence classes.
\mn
\ub{Case 2}:  $\mu > \aleph_0$.

Choose $<_\alpha$ be a well ordering of ${}^\alpha 2$ and let $E$ be the
following equivalence relation on $[{}^\mu 2]^{n^*}$: if $w_0 =
\{\eta_\ell:\ell < n^*\},w_2 = \{\nu_\ell:\ell < n^*\}$ are members of
$[{}^\mu 2]^{n^*}$ both listed in lexicographic increasing order then:
$w_1 E w_2$ iff for any $\ell_1 < \ell_2 < n^*$ and $\ell_3 < \ell_4 < n^*$
we have
\mr
\item "{$(a)$}"  $\ell g(\eta_{\ell_1} \cap \eta_{\ell_2}) 
< \ell g(\eta_{\ell_3} \cap \eta_{\ell_4}) \Leftrightarrow \ell g
(\nu_{\ell_1} \cap \nu_{\ell_2}) < \ell g(\nu_{\ell_3} \cap \nu_{\ell_4})$
\sn
\item "{$(b)$}"  $\eta_{\ell_3} \restriction \ell g(\eta_{\ell_1} \cap 
\eta_{\ell_2}) <_{\ell g(\eta_{\ell_1} \cap
\eta_{\ell_2})} \eta_{\ell_4} \restriction \ell g(\eta_{\ell_1} \cap
\eta_{\ell_2}) \Leftrightarrow \nu_{\ell_3} \restriction \ell g(\nu_{\ell_1} 
\cap \nu_{\ell_2}) <_{\ell g(\nu_{\ell_1} \cap \nu_{\ell_2})}
\nu_{\ell_4} \restriction \ell g(\nu_{\ell_1} \cap \nu_{\ell_1})$.
\endroster
\endremark
\newpage

\shlhetal

\newpage
    
REFERENCES.  
\bibliographystyle{lit-plain}
\bibliography{lista,listb,listx,listf,liste}

\enddocument

\bye